\documentclass[11pt,  oneside]{article}
\usepackage[a4paper, margin={1.2in}]{geometry}

\usepackage[utf8]{inputenc}
\usepackage[T1]{fontenc}
\usepackage{tgpagella}
\usepackage{amssymb}
\usepackage{amsfonts}
\usepackage{amsmath}
\usepackage{amsthm}
\usepackage{mathrsfs}
\usepackage[table]{xcolor}
\usepackage{booktabs}
\usepackage{multirow}
\usepackage{color}
\usepackage{amscd}
\usepackage{bussproofs}
\usepackage[all,color]{xy}
\usepackage{comment}

\newtheorem{tw}{Theorem} 
\newtheorem{wniosek}[tw]{Corollary} 
\newtheorem{lemat}[tw]{Lemma} 
\newtheorem{fakt}[tw]{Fact} 
\newtheorem{stwierdzenie}[tw]{Proposition} 

\theoremstyle{definition}
\newtheorem{definicja}[tw]{Definition} 
\newtheorem{konwencja}[tw]{Convention} 

\newtheorem{przyklad}[tw]{Example} 
\newtheorem{uwaga}[tw]{Remark} 

\newcommand{\qcr}[1]{\ulcorner #1 \urcorner}

\newcommand{\CT}{\textnormal{CT}}
\newcommand{\PA}{\textnormal{PA}}
\newcommand{\Sent}{\textnormal{Sent}}

\newcommand{\set}[2]{\{#1 \ \ | \ \ #2 \}}
\newcommand{\pair}[2]{\langle #1, #2\rangle}

\newcommand{\form}{\textnormal{Form}}
\newcommand{\CS}{\textnormal{CS}}

\newcommand{\ElDiag}{\textnormal{ElDiag}}
\newcommand{\Th}{\textnormal{Th}}
\newcommand{\Form}{\form}

\newcommand{\was}[1]{\{#1\}}

\newcommand{\EA}{\textnormal{EA}}

\newcommand{\UTB}{\textnormal{UTB}}
\newcommand{\ClTerm}{\textnormal{ClTerm}}
\newcommand{\Ass}{\textnormal{Asn}}

\newcommand{\restr}[1]{{\upharpoonright_{#1}}}

\newcommand{\Sat}{\textnormal{Sat}}
\newcommand{\True}{\textnormal{True}}
\newcommand{\Con}{\textnormal{Con}}
\newcommand{\Tr}{\textnormal{Tr}}

\newcommand{\Proof}{\textnormal{Proof}}
\newcommand{\dpt}{\textnormal{dp}}

\newcommand{\Lang}{\mathcal{L}}
\newcommand{\Ind}{\textnormal{Ind}}

\newcommand{\M}{\mathcal{M}}

\newcommand{\N}{\mathcal{N}}
\newcommand{\K}{\mathcal{K}}

\newcommand{\tuple}[1]{\langle #1 \rangle}

\newcommand{\Def}{\textnormal{Def}}

\newcommand{\num}[1]{\underline{#1}}

\newcommand{\DEF}{\textnormal{DEF}}
\newcommand{\UDEF}{\textnormal{UDEF}}
\newcommand{\Diag}{\textnormal{Diag}}
\newcommand{\USB}{\textnormal{USB}}

\newcommand{\pdpt}{\textnormal{pdp}}
\newcommand{\SAT}{\textnormal{SAT}}
\newcommand{\RSAT}{\textnormal{RSAT}}

\newcommand{\Aut}{\textnormal{Aut}}

\newcommand{\val}[1]{{#1}^{\circ}}

\newcommand{\TB}{\textnormal{TB}}

\newcommand{\rem}[1]{#1}

\mathchardef\mhyphen="2D

\title{Universal properties of truth}
\author{Mateusz Łełyk, Bartosz Wcisło}

\begin{document}
	
	\maketitle
	
	\section{Introduction}
	
	Many types of arguments in mathematical logic (especially in foundations) make an important use of a truth predicate of some form. For example, when one wants to prove a consistency of one theory in its extension, one often proceeds by defining a truth predicate for a class of formulae and then inductively shows that each sentence from this class which is provable in the weaker theory is true. Similarly, if one wants to deduce finite axiomatizability of an infinite theory, one often makes use of a definable truth predicate. Last but not least, the existence of a truth predicate for some language $\Lang$ is often used to show that a given model of an $\Lang$-theory is $\Lang$-recursively saturated.
	
	\rem{In most general terms, the current paper investigates whether the use of the truth predicate in such arguments is \textit{essential}}.
	More concretely, we focus on model-theoretic properties and ask whether for a given property $P$, the existence of a definable truth predicate is necessary to uniformly impose this property. One of the clearest illustrations here is provided by the property of \textit{imposing recursive saturation}: it is a very basic fact that every nonstandard model $\M\models \PA$ which has a partial inductive nonstandard satisfaction class (i.e. a satisfaction predicate for all formulae of some nonstandard complexity which satisfy full induction) is recursively saturated (a very simple overspill argument is presented in \cite{kaye}). Since the theory of partial inductive nonstandard satisfaction classes is axiomatizable, we have an example of a theory, call it $\UTB$, which \textit{imposes $\Lang_{\PA}$-recursive saturation}. In general, we say that a theory $U$ imposes $\Lang$-recursive saturation iff the $\Lang$-reduct of an arbitrary nonstandard model of $U$ is recursively saturated. Quite surprisingly there is a good sense in which the use of the $\UTB$-truth predicate is essential in this context. More precisely, we have the following theorem due to Roman Kossak \cite[Theorem 2.4]{kossak_cztery_problematy}:
	
	\begin{tw}\label{tw_kossak}
		Suppose that $U$ is a theory in a language $\Lang$ which extends $\PA$ and proves all instantiations of the induction scheme with $\Lang$-formulae. If $U$ imposes $\Lang_{\PA}$-recursive saturation, then in every $\M\models U$ there exists a definable partial satisfaction class.
	\end{tw}  
	
	Thus Kossak's theorem shows that, for a reasonable class of theories, having a definable satisfaction predicate is necessary to impose recursive saturation.
	
	Another example of a property that seems even more directly connected with the definability of the notion of truth is \textit{imposing elementary equivalence}. Consider the weak truth theory $\TB^-$: it is obtained by  adding $\PA$ all the axioms of the form $T(\qcr{\phi}) \equiv \phi$ (where $T$ is a fresh predicate). Then $\TB^-$ has the following property: whenever $\M \models \TB^-$, $\N\models \TB^-$ and $\M$ is a submodel of $\N$, then $\M$ and $\N$ satisfy the same sentences of $\Lang_{\PA}$. This can be abstracted into a property of a general theory: we say that a theory $U$ imposes $\Lang$-elementary equivalence iff whenever $\M,\N\models U$ and $\M\subseteq \N$, then $\M$ and $\N$ satisfy the same $\Lang$-sentences. Then we can ask whether having a definable truth predicate for $\Lang$-sentences is necessary for a (sufficiently strong) theory $U$ to impose $\Lang$-elementary equivalence. It turns out that $\TB^-$ is indeed essential to ensure this property (Theorem \ref{tw_impos_elem}, point 3).

	When we start investigating the main question, one quickly encounters a number of semantic properties imposed by certain natural truth theories over $\PA$ for which those truth predicates are not essential. However there seem to be quite natural candidates for axiomatic theories of truth-related notions which do provide such exact characterisations. One example of such a phenomenon is the notion of definability: whenever $\M, \N$ are models of $\UTB^-$  (the theory of the uniform Tarski biconditionals without induction) and $\M$ is a submodel of $\N$, they have the same arithmetically definable elements. This can be again defined as an abstract property of a theory. However, in order for that property to hold it is not essential that we have a truth predicate to our disposal. Nevertheless, we can introduce an axiomatic theory of definability, $\DEF^-$, which indeed corresponds to that semantic property. So we can expand our initial question as follows: given a truth-related semantic property, find an axiomatic theory whose use is essential in ensuring that property. 
	
	As we already mentioned, imposing definability is in general not enough to define a predicate satisfying $\UTB^-$. This can be in fact rephrased by saying that in $\DEF^-$ one cannot define a predicate satisfying $\UTB^-$ (while preserving the arithmetical part of the language). Therefore, it becomes natural to ask about definability relations between theories corresponding to truth-like semantic properties. 
	
	The definability notions which we analyse, come in many flavours: we can consider them with or without parametres. Similarly, we can require that one theory defines another over every model, possibly with different formulae or we can require the definability to hold uniformly in the theory. Finally, we can require that the definability is non-uniform between models, but uniformly bounded in complexity. Those distinctions apply to both types of questions we discussed above: we can consider various definability notions of a theory corresponding to a semantic property or analyse distinct reductions between those theories. 
	
	Finally, it sometimes happen that the answer to an analysed question depends on the base theories. In a similar manner, the choice of the base theory may vastly simplify the arguments. In general, we will try to work in a general framework of arbitrary sequential theories, but in some cases we can get a clearer picture or more elegant proofs, when we accept some additional natural assumptions on the base theories (notably, we can choose $\PA$). This adds another dimension to our investigations.
	
The plan of the paper is as follows.
	In section 2 we introduce all the relevant notions, including the axiomatic theories of truth ($\TB$), \rem{definability ($\DEF$)} and satisfaction ($\USB$ or $\UTB$) in the general context of an arbitrary sequential theory. The theory $\DEF$ is introduced in this paper for the first time. Moreover, we state the definitions of two definability relations between theories, which we call syntactical and semantical definability. Finally we introduce abstract model-theoretical properties of general theories that will be used to characterize axiomatic theories of semantical notions. Each of the properties is introduced in two variants: uniform and non-uniform. 
	
	In section 3 each of the axiomatic theories is paired with the appropriate model-theoretical notion: we show that, up to syntactical definability, axiomatic theories of truth, definability
	and satisfaction can be characterized via the notions of imposing elementary equivalence, preserving definability and imposing elementarity, respectively (Theorem \ref{tw_impos_elem}). The theorems account for coordinate-free characterizations of these axiomatic theories, in the sense of \cite{VisserR}. In Theorems \ref{tw_nonuniform_def_tb} and \ref{tw_nonuniform_impos_elem_utb} we prove that the above characterization holds also for semantical definability and the non-uniform analogues of the model-theoretical properties.  
	
	Section 4 explores definability relations between various axiomatic theories of semantic notions. These results can be also viewed as  generalisations of some previously known non-definability results between the axiomatic theories of truth and satisfaction (see \cite[Corollary 3.8]{wcislyk_models_weak} and \cite[Theorem 5.9]{viss19enayat}). 
	We show that satisfaction is not definable from truth over any sequential base theory (Theorem \ref{tw_nondef_sat_in_def_tb}). This solves problem 5.10 in \cite{viss19enayat}. Then this picture is enriched with an axiomatic theory of definability ($\DEF$): we show that, except for some degenerate cases\footnote{We classify as \rem{"degenerate"} the examples of sequential theories which, roughly speaking, are inconsistent with the theory of extensional adjunctive sets.}, the theory situates itself strictly between $\TB$ and $\USB$ (Theorems \ref{tw_nondef_sat_in_def_tb} and \ref{tw_nondef_def_in_tb}). It is shown that each of the separations ($\DEF$ from $\TB$ and $\USB$ from $\DEF$) hold in the strict sense: each (non-degenerate) sequential theory has a model which can be expanded to a model of $\TB$ ($\DEF$) such that in the expanded model $\DEF$ ($\USB$, resp.) is not definable. The proof for $\TB$ vs. $\DEF$ case uses a construction of a model in which the definability hierarchy does not collapse. The construction works uniformly for any (non-degenerate) sequential theory $U$, and to our best knowledge this is an original result in the model theory of sequential theories (Lemma \ref{lem_hierarchical_model}).  As a corollary we conclude that, unless $U$ is degenerate, $\DEF$ does not have a restricted axiomatization modulo $U$. A similar argument shows that for every sequential $U$, $\TB$ does not have a restricted axiomatization modulo $U$, which solves a problem posed in \cite[Open Question 5.7]{viss19enayat}. In the final two subsections we study the definability relations between the above mentioned theories and two axiomatic theories of Skolem functions (uniform and non-uniform; Subsection \ref{subsect_skolem}) and the special
	 situation with $\PA$ as the base theory (Subsection \ref{subsect_PA}). In the former, we use Kaye-Kossak style results (from \cite{kossak_1990}, \cite{kossak_cor_to_1990} and \cite{kaye_model_char_PA}) about models of fragments of $\PA$ to construct models of $\USB$ in which there is no definable witness-picking function for arithmetical formulae (Theorems \ref{stw_nondef_H_usb} and \ref{tw_nondef_US_in_UTB}). In the latter, we show that various structural results about models of $\PA$ provide simpler and stronger (Subsection \ref{subsect_PA_def_inno}) separability arguments for $\TB$, $\DEF$ and $\USB$.
	
	In section 5 we study the property of imposing recursive saturation and  possible strengthenings of Kossak's theorem (Theorem \ref{tw_kossak}). Using a recursion-theoretic argument (involving an arithmetisation of a construction of an $\omega$-generic set intertwined with the Jump Inversion Theorem) we show that there is a theory in a finite language which imposes recursive saturation but does not syntactically define $\UTB$ (Theorem \ref{th_theory_finlang_with_rec_sat_without_utb}). To our best knowledge this is the first example of a theory in a finite language which defines $\UTB$ semantically but not syntactically. Secondly, we show that Kossak's Theorem does generalise to theories without full induction (but extending $\PA$): we prove that each sequential theory in a countable language which extends $\PA$ and imposes $\Lang_{\PA}$-recursive saturation semantically defines $\UTB$ for $\Lang_{\PA}$ (Theorem \ref{th_recursive_saturation_implies_truth}).\footnote{Here in $\UTB$ we do not assume the induction scheme for formulae with the satisfaction predicate.} The proof is very different from the original Kossak's reasoning and uses iterated omitting types theorem in the spirit of \cite{casanovas_farre}. As a corollary we define a theory $\RSAT$ which is universal for imposing $\Lang_{\PA}$-recursive saturation in the following sense: a sequential theory $U$ in a countable language and extending $\PA$ imposes $\Lang_{\PA}$-recursive saturation if and only if $U$ semantically defines $\RSAT$.
	
	\rem{Section 6 contains a list of open problems. Most of them are connected either to some questions in model theory of fragments of $\PA$ or interpretability and definability in finitely axiomatized theories.}

	\section{Preliminaries}\label{sect_prelim}
	
	\subsection{Model-theoretical properties}
	
	For general first-order theories the arithmetical $\Sigma_n\setminus\Pi_n$ hierarchy need not make much sense. In what follows we shall employ the quantifier alternation depth hierarchy, introduced in \cite{Buss2015} and in a similar context used in \cite{viss_sivs}.
	
	\begin{definicja}[Depth of quantifier alternation]\label{def_complexity}
		Let $\Lang$ be any first-order language and let $\textsc{AT}$ denote the set of atomic formulae of $\Lang$. $\Sigma_n^*(\Lang)$ and $\Pi_n^*(\Lang)$ are defined by simultaneous induction as follows  (below, and whenever it is clear from the context, we omit the reference to $\Lang$):
		\begin{align*}
			\Sigma^*_0 &= \Pi_0^* := \emptyset\\
			\Sigma_{n+1}^*&:= \textsc{AT}\ \ | \ \ \exists v \Sigma_{n+1}^* \ \ | \ \ \Sigma^*_{n+1}\wedge \Sigma^*_{n+1} \ \ | \ \ \Sigma^*_{n+1}\vee \Sigma^*_{n+1} \ \ | \ \ \neg\Pi^*_{n+1} \ \ | \forall v \Pi^*_n | \\
			\Pi_{n+1}^*&:= \textsc{AT}\ \ | \ \ \forall v \Pi_{n+1}^* \ \ | \ \ \Pi^*_{n+1}\wedge \Pi^*_{n+1} \ \ | \ \ \Pi^*_{n+1}\vee \Pi^*_{n+1} \ \ | \ \ \neg\Sigma^*_{n+1} \ \ | \exists v \Sigma^*_n | \\
        \Delta_{n+1}^*&:= \Sigma_{n+1}^*\cap \Pi_{n+1}^*
		\end{align*}
		For a formula $\phi$, $\Sigma^*(\phi)$ $(\Pi^*(\phi), \Delta^*(\phi))$  denotes the $\Sigma^*_n$ ($\Pi^*_n, \Delta_n^*$) complexity of $\phi$.

        We will need two more fine-grained measures of complexity of formulae. We call the first one \emph{depth} ($\dpt$) and the second \emph{pure depth} ($\pdpt$). In short, the depth of $\phi$ is the maximal number of connectives and quantifiers on a branch in the syntactic tree of $\phi$. The pure depth of $\phi$ is the maximal length of a branch in the syntactic tree of $\phi$, where terms occuring in $\phi$ are unravelled. Formally, we define $\dpt$ and $\pdpt$ recursively:
        \begin{align*}
        \dpt(\phi) &= 0, \textnormal{ if } \phi\in \textsc{AT}\\
        \dpt(\neg\phi)&= \dpt(\exists x \phi) = \dpt(\forall x \phi) = \dpt(\phi)+1\\
        \dpt(\phi\wedge\psi) &= \dpt(\phi\vee\psi) = \max\{\dpt(\phi),\dpt(\psi)\}+1 
        \end{align*}
        Pure depth is defined first for terms: $\pdpt(t)=0$, if  $t$ is a constant or a variable and $\pdpt(f(t_1,\ldots,t_n)) = \max\{\pdpt(t_1),\ldots,\pdpt(t_n)\}+1.$ Then $\pdpt(\phi)$ is defined as $\dpt(\phi)$, except for the case of atomic formulae, where we put $\pdpt(R(t_1,\ldots,t_n)) = \max\{\pdpt(t_1),\ldots,\pdpt(t_n)\}$.
	\end{definicja}
	
	In the paper we only consider first-order languages. If $\M$ is any model of a language $\Lang$, $A\subseteq M$, then $\Lang_A$ denotes the extension of $\Lang$ with constants naming every element of $A$. If $U$ is any first order theory, then $\Lang_U$ denotes its language. If $\M$ is any model and $\Lang$ any language, then $\M\restr{\Lang}$ denotes the $\Lang$-reduct of $\M$.

	\begin{definicja}
		Let $\Lang$ an arbitrary language.
		\begin{enumerate}
			\item $n\mhyphen \Diag(\M)$ denotes the set of $\Sigma^*_n(\Lang_M)$ sentences, that are true in $\M$. We put $\ElDiag(\M) = \bigcup_{n\in\mathbb{N}}n\mhyphen \Diag(\M)$. We write $\M\preceq_{n} \N$ as an abbreviation of $n\mhyphen\Diag(\M)\subseteq n\mhyphen\Diag(\N)$ and $\M\preceq \N$ for $\ElDiag(\M)\subseteq \ElDiag(\N)$.
			\item $n\mhyphen \Th(\M)$ denotes the set of $\Sigma^*_n(\Lang)$ that are true in $\M$.  We put $\Th(\M) = \bigcup_{n\in\mathbb{N}}n\mhyphen \Th(\M)$. We write $\M\equiv_{n} \N$ as an abbreviation of $n\mhyphen\Th(\M)= n\mhyphen\Th(\N)$ and $\M\equiv \N$ for $\Th(\M)= \Th(\N)$.
			\item $n\mhyphen \Def(\M)$ denotes the set of $\Sigma^*_n(\Lang)$ definable elements of $\M$. We put $\Def(\M) = \bigcup_n n\mhyphen\Def(\M)$.
			\item  $n\mhyphen \DEF(\M)$ denotes the set $\set{\pair{\theta(x)}{a}}{\M\models \theta(a)\wedge \exists!x\theta(x)}$.  We put $\DEF(\M) = \bigcup_n n\mhyphen\DEF(\M)$.
		\end{enumerate}   
	\end{definicja}
	
	\begin{definicja}
		Let $U$ be any theory and $\Lang$ be an arbitrary language.
		\begin{enumerate}
			\item We say that $U$ \emph{imposes $\Lang$-elementarity} (\emph{$\Lang$-elementary equivalence}) if for every $\M\models U$ there exists $n\in\mathbb{N}$ such that for all $\N\models U$, 
			\begin{align}\label{impos_elem_def}\tag{E$^{\M, \N}_n$}
				\M\preceq_{n}\N &\Rightarrow \M\restr{\Lang}\preceq\N\restr{\Lang}\\
				(\M\equiv_{n}\N&\Rightarrow \M\restr{\Lang} \equiv \N\restr{\Lang}). \label{impos_elem_equiv_def}\tag{EE$^{\M,\N}_n$}
			\end{align}
			
			\item We say that $U$ \textit{uniformly imposes $\Lang$-elementarity} (\textit{$\Lang$-elementary equivalence}) if there is an $n$ such that for every $\M,\N\models U$, \eqref{impos_elem_def} holds (\eqref{impos_elem_equiv_def} holds, respectively).
			
			\rem{	\item We say that $U$ \textit{imposes equality of $\Lang$-definables} (\textit{preserves $\Lang$-definability}) if for every $\M\models U$ there is $n\in\mathbb{N}$ such that for every $\N\models U$
				\begin{align}
					n\mhyphen\Def(\M)= n\mhyphen\Def(\N) &\Rightarrow \Def(\M\restr{\Lang}) = \Def(\N\restr{\Lang}).\label{impos_equal_def}\tag{ED$^{\M,\N}_n$}\\
					\bigl(n\mhyphen \DEF(\M)= n\mhyphen\DEF(\N) &\Rightarrow \DEF(\M\restr{\Lang}) = \DEF(\N\restr{\Lang})\bigr) \label{preserv_def}\tag{PD$^{\M,\N}_n$}
				\end{align}
			}	
			\item We say that $U$ \textit{uniformly imposes equality of $\Lang$-definables} (\textit{preserves $\Lang$-definability}) if there is an $n$ such that for every $\M,\N\models U$, \eqref{impos_equal_def} holds (\eqref{preserv_def} holds, respectively).
			
			\item We say that $U$ \textit{imposes recursive $\Lang$-saturation} if for every $\M\models U$, $\M\restr{\Lang}$ is recursively saturated.
		\end{enumerate}
	\end{definicja}
	
	\subsection{Various grades of definability}
	
	In what follows the definitions of translation and interpretation are borrowed from \cite{viss_sivs}. We use $K: U\lhd V$ to denote the fact that $K$ is a relative interpretation of $U$ in $V$. $id_U$ denotes the identity interpretation $U\lhd U$ and we omit the reference to $U$ if it is clear from context. Let $\Theta_0$, $\Theta_1$ be two signatures and $\Lang$ a language. We say that the translation $\tau:\Theta_0 \rightarrow \Theta_1$ is $\Lang$-conservative if for every formula $\phi\in\Lang$, $\phi^{\tau}=\phi$ and we say that $\tau$ is direct if it is one-dimensional, unrelativized and it translates identity to identity.
	
	\begin{definicja}
		We say that $U$ \emph{syntactically defines} $V$ \textit{modulo} $\Lang$, $V\leq_{\Lang} U$, if $U$ interprets $V$ via a direct, parameter-free and $\Lang$-conservative interpretation.
		
		We say that $U$ \rem{\emph{semantically defines} $V$ \textit{modulo} $\Lang$, $V\leq_{\Lang}^m U$} if for every $\M\models U$, $\ElDiag(\M)$ syntactically defines $V$ modulo $\Lang_{M}$.

		We say that $U$ \emph{semantically parameter-free defines $V$ modulo} $\Lang$ if every consistent and complete extension of $U$ syntactically defines $V$ modulo $\Lang$. 
	\end{definicja}
	
	All the definability results in this paper consider the situation in which the defined-modulo-$\Lang$ theory is formulated in a language extending $\Lang$ with a single (at most binary) predicate $R$. In such a context it is easy to see that $U$ \textit{syntactically defines} $V$ \textit{modulo $\Lang$} if and only if there is a formula $\phi_R(\bar{x})\in\Lang_U$ such that for every $\Phi$-axiom of $V$
	\[U\vdash \Phi[\phi_R(\bar{t})/R(\bar{t})],\]
	where $\Phi[\phi_R(\bar{t})/R(\bar{t})]$ denotes the result of substituting $\phi_R(\bar{t})$ for each occurrence of $R(\bar{t})$, perhaps preceded by renaming of the bounded variables in $\phi_R(\bar{t})$ so as to avoid unwanted variable capture. 
	
	\begin{uwaga}
		It is very easy to observe that if $U$ syntactically defines $V$, then $U$ semantically defines $V$ and that this reverses if $V$ is finite. However, in a general context, it is possible for $U$ to semantically, but not syntactically define $V$. In such a scenario, in every model $\M$ of $U$ one can define intepretations for the predicates in $\Lang_V$ so as to turn $\M$ into a model of $V$, but there is no single fixed interpretation that works across all models. Our paper provides many examples of this phenomenon in the realm of axiomatic theories of semantical notions.
	\end{uwaga}

	\subsection{Axiomatic theories of semantical notions}\label{sect_prelim_truth}
	
	An axiomatic theory of truth (satisfaction, definability,...) splits into two parts: the base theory and some axioms specific for the notion of truth. An important restriction is that the syntax is encoded uniquely by the base theory. For the very basic theories of truth studied in this paper to make sense already very moderate syntax theory is sufficient: as shown in \cite{viss19enayat}, theories of truth make sense already for theories which are able to interpret the basic theory of successor and theories of satisfaction make sense over Vaught theories. We shall start from a stronger assumption that the base theory is sequential, which already covers a wide range of important first-order theories.   
	
	\subsubsection{Sequential theories}
	
	A theory is sequential if it admits a reasonably well-behaved notion of a finite sequence. More precisely: a sequential theory needs to prove that there exists an empty sequence and that every sequence $s$ can be extended by appending a given element from the universe at the end of $s$. Here comes the formal definition:
	
	\begin{definicja}[Adjunctive sets, sequentiality]\label{defi_sekwencyjność}
		The theory $\textsc{AS}$ is formulated in the language with one binary relational symbol $\in$ and has as axioms
		\begin{itemize}
			\item[AS1] $\exists x \forall y \bigl(y\notin x\bigr)$.
			\item[AS2] $\forall x \forall y \exists z \forall w \bigl(w\in z \equiv w\in x \vee w = y\bigr)$.
		\end{itemize}
		A theory $U$ is \emph{sequential} if there is a direct interpretation of $\textsc{AS}$ in $U$. The translation of $\in$ under this interpretation will be denoted $\in_U$. A model $\M$ is \emph{sequential} if $\ElDiag(\M)$ is sequential. We write $\in_{\M}$ instead of $\in_{\ElDiag(\M)}$. If $\M$ is sequential and $c\in M$, then $c^{\M}$ denotes the set $\set{a\in M}{\M\models a\in_{\M} c}$.
		
		For later purposes let us introduce also the theory of \textit{adjunctive sets with extensionality}, ASE, which is an extension of AS with the axiom $\forall x\forall y \bigl(\forall z (z\in x\equiv z\in y)\rightarrow x=y\bigr)$.
	\end{definicja}
	
	In this paper we shall restrict ourselves to sequential theories which directly interpret $\textsc{AS}$ via a parameter-free interpretation.
	We leave the problem of adding parameters for further investigations.
	
	Let us recall some prominent sequential theories: $\PA$ denotes Peano Arithmetic which we take to be formulated in the language $\Lang_{\PA}$ of ordered rings $\was{0,1,+,\times, <}$. I$\Sigma_n$ denotes the fragment of $\PA$ consisting of axioms of induction uniquely for formulae in $\Sigma_n$ form (according to the standard definition in which $\Sigma_0$ ($\Delta_0$) consists of all formulae in which all quantifiers are bounded) and B$\Sigma_n$ extends I$\Delta_0$ with collection axioms for $\Sigma_n$ formulae. Abusing the notation, we use $\EA$ to denote the theory in the language $\Lang_{\PA}$ extending $I\Delta_0$ with the axiom $\exp$ saying that the exponential function is total. All the details and basic facts about these theories can be found in \cite{kaye} and \cite{hapu98}. Finally, $S^1_2$ is the Buss's weak arithmetic of P-Time computability (consult \cite[Chapter V.4.b, Definition 4.4]{hapu98} for details).
	
	The notion of sequentiality is extensively discussed e.g. in \cite{viss_sivs}. Two most important for the current project features of sequential theories are summarized in the following two facts (Fact \ref{sequentiality_arithmetic} and Fact \ref{fakt_partial_truth}). 
	
	\begin{fakt}\label{sequentiality_arithmetic}
		If $U$ is sequential, then $U$ relatively interprets $S^1_2$ via a one-dimensional, parameter-free interpretation (but we allow the equality to be re-defined).
	\end{fakt} 
	
	To see how this is possible, observe that in AS one can first interpret ASE (by introducing new identity relation which glues together sets with the same elements). Then inside this interpretation one can define von Neumann ordinal numbers in the standard way. By restricting the domain to those ordinals on which multiplication and addition behaves nicely, one obtains a relative interpretation of Robinson's arithmetic $Q$. The fact that $Q$ interprets $S^1_2$ is well-known (see e.g. \cite{FerrQ}). 
	
	\begin{konwencja}\label{konwencja_arytmetyka_w sekwencji}
		We make use of the following abbreviations:
		\begin{itemize}
			\item We assume fixed G\"odel coding of languages in consideration. For a fixed formula $\phi$, $\qcr{\phi}$ denotes the G\"odel code of $\phi$.
			\item $\num{n}$ denotes the binary numeral naming $n$, to wit, the expression:
			\[\num{a_0} + \num{2}\cdot (\num{a_1} + \num{2}\cdot (\num{a_2} + \ldots(\num{a_{k-1}} + \num{2}\cdot \num{a_{k}})\ldots)),\] 
			where $a_i\in\was{0,1}$, $\num{0} = 0$, $\num{1} = 1$ and $\num{2} = (1+1)$ and $(a_i)_{i\leq k}$ is a binary expansion of $n$.
			\item Throughout the paper, when working in a sequential theory $U$ with a fixed interpretation $N$ of $S^1_2$, $E$ shall always denote the definable  equivalence relation  which is used by $N$ to translate the identity relation from the arithmetical signature. For any $n\in\omega$, we  write $n(x)$ for the $N$-translation of the formula $x = \num{n}$. We stress that, in the context of a general sequential theory, $\phi(\num{n})$ should be understood contextually as $\exists x \bigl(n(x) \wedge \phi(x)\bigr)$, where $n(x)$ is a predicate expressing the property of "being the $n$-th ordinal". We stress that $n(x)$ need not define any particular object, although it uniquely determines an $E$-equivalence class. Putting together this and the previous conventions, $\num{\qcr{\phi}}$ will denote a predicate corresponding to the G\"odel code of a formula $\phi$.
			\item Similarly, in the context of a sequential theory $U$ and $N: S^1_2\lhd U$, we shall use arithmetical predicates such as the ordering ($<$), "$x$ is a sentence of a language $\Lang$" ($\Sent_{\Lang}$), "$x$ is a formula of a language $\Lang$" ($\Form_{\Lang}$), "$x$ is a proof of $y$ in pure first order logic" ($\Proof_{\emptyset}(x,y)$), "$x$ is a $\Sigma^*_n$ complexity of $\phi$" ($x=\Sigma_n^*(\phi)$), ''$\alpha$ is an assignment'' ($\alpha\in\Ass$) having in mind their $N$-translations to the language $\Lang_U$. We also assume that these predicates apply uniquely to members of the domain of $N$.
			\item To keep the good balance between precision and readability, we often abbreviate the use of a syntactic operation by using $\qcr{\cdot}$ (stretching the notation introduced above). For example if $\phi$ is understood to be a sentence, then $\qcr{\exists x \phi}$ will denote the result of prefixing $\phi$ with an existential quantifier and then variable $x$. 	 For substitutions of binary numeral we use the dot notation: $\qcr{\phi(\dot{x}/v)}$ denotes the result of substituting the binary numeral naming $x$ for the variable $v$ in $\phi$. We stress that we often treat $\qcr{\phi(\dot{x}/v)}$ as a function in three variables, $x$, $v$ and $\phi$ and skip the reference to $v$ if it is clear from context.
		\end{itemize}	
		
	\end{konwencja}

	\begin{fakt}\label{fakt_partial_truth}
		If $U$ is a sequential theory, then for each $n$ there is a formula $\Sat_n(x,y)$ which, provably in $U$ satisfies the Tarskian compositional clauses for formulae of complexity $\Delta^*_n$ whose pure depth belong to a certain $U$-provable cut $J$.
	\end{fakt}
    More precisely, for a fixed sequential theory $U$ in a finite language $\Lang$, $\CS^-(\Lang)(\phi)$ denotes the following formula of $\Lang\cup S$, where $S$ is a fresh binary predicate (the initial disjunction $\bigvee_{R\in \Lang}$ and $\bigvee_{f\in\Lang}$ range over all relational symbols and all function symbols from $\Lang$, respectively (we treat constants as $0$-ary function symbols)):
    \begin{align*}
    \forall \alpha\in\Ass\bigl[S(\phi,\alpha) \equiv& \bigl[\bigvee_{R\in \Lang}\exists v_0,\ldots, v_n \bigl(\phi = \qcr{R(v_0,\ldots,v_n)}\wedge R(\alpha(v_0),\ldots,\alpha(v_n)\bigr)  \\
    &\vee \bigvee_{f\in \Lang}\exists v_0,\ldots, v_n \bigl(\phi = \qcr{v_{0} = f(v_1,\ldots,v_n)}\wedge \alpha(v_0) = f(\alpha(v_1),\ldots,\alpha(v_n))\bigr)\\
    &\vee \exists \phi_0,\phi_1\bigl(\phi = \qcr{\phi_0\wedge\phi_1}\wedge (S(\phi_0,\alpha)\wedge S(\phi_1,\alpha))\bigr)\\
    &\vee \exists \phi_0,\phi_1\bigl(\phi = \qcr{\phi_0\vee\phi_1}\wedge (S(\phi_0,\alpha)\vee S(\phi_1,\alpha))\bigr)\\
    &\vee \exists \phi_0 \bigl(\phi = \qcr{\neg\phi_0}\wedge \neg S(\phi_0,\alpha)\bigr)\\
    &\vee \exists \phi_0\exists v_0\bigl(\phi = \qcr{\forall v_0}\wedge \forall\beta \bigl(\beta\sim_{v_0}\alpha \rightarrow S(\phi_0,\beta)\bigr)\\
    &\vee \exists \phi_0\exists v_0\bigl(\phi = \qcr{\exists v_0}\wedge \exists\beta \bigl(\beta\sim_{v_0}\alpha \wedge S(\phi_0,\beta)\bigr)\bigr]\bigr].
    \end{align*}
	The results of \cite{viss_sivs}[Section 2] show that for every $n$ there exists a $U$-provable cut $J_n$ and a formula $\Sat_n(x,y)$ such that
   \[U\vdash \forall \phi\bigl(J_n(\pdpt(\phi))\wedge \Delta^*_n(\phi)\rightarrow \CS^-(\Lang)(\phi^{te})[\Sat_n(x,y)/S(x,y)]\bigr),\]
	where $\phi^{te}$ is the canonical formula equivalent to $\phi$, in which term symbols occur uniquely in the context $y=f(x)$, where $y$ and $x$ are variables. In what follows, $\SAT^{x,y}_n[\psi]$ denotes the formula 
 \[\forall \phi\bigl(J_n(\pdpt(\phi))\wedge \Delta^*_n(\phi)\rightarrow \CS^-(\Lang)(\phi^{te})[\psi(x,y)/S(x,y)]\bigr)\]
	We note that $\SAT_n^{x,y}[\psi]$ might contain free variables, if $\psi$ has some free variables other that $x$ and $y$. The last definition generalizes the notion of a coded set from the context of fragments of $\PA$ to arbitrary sequential theories.
	\begin{definicja}\label{defi_kodowanie}
		Let $\M$ be a sequential model and $K:S^1_2\lhd \ElDiag(\M)$. We say that a set $A\subseteq \omega$ is \emph{coded} in $\M$ if for some $c\in M$ and every $n\in\omega$, $\M\models \exists x \bigl(n(x)\wedge x\in_{\M} c\bigr)$ if and only if $n\in A$.
	\end{definicja}

	\subsubsection{Truth, satisfaction, definability}
	
	In this section we have a sequential theory $U$ in the background, together with the distinguished one-dimensional parameter-free interpretation $K:S^1_2\lhd U$. We assume that the language of $U$ is r.e. and the codes of formulae are given by this interpretation. Hence $\qcr{\phi}$ denotes the natural number coding $\phi$ and $\num{\qcr{\phi}}$ is the predicate picking the $E$-equivalence class corresponding to $\qcr{\phi}$:
	
	\paragraph{Truth} $\TB_K(\Lang)$ extends $U$ with all axioms of the form
	\[T(\num{\qcr{\phi}})\equiv \phi,\]
	where $\phi$ is a $\Lang$-sentence and $T$ is a fresh unary predicate.	
	
	\paragraph{Definability}
	
	$\DEF_K(\Lang)$ extends $U$ with all axioms of the form:
	\[\forall y \bigl(D(\num{\qcr{\phi(x)}}, y)\equiv \exists! x \phi(x)\wedge \phi(y)\bigr),\]
	for $\phi\in\Lang$ and $D$- a fresh binary predicate.
	
	\paragraph{Satisfaction}  $\USB_K(\Lang)$ extends $U$ with all axioms of the form
	\[\forall x \bigl(S(\num{\qcr{\phi}}, x)\equiv \phi(x)\bigr),\]
	where $\phi$ is a $\Lang$-formula and $S$ is a fresh binary predicate.
	
	\begin{uwaga}
		The reference to $U$ in $\TB_K$, $\DEF_K$ and $\USB_K$ is hidden in the interpretation $K$: as in \cite{viss_sivs} we take both the interpreting and the interpreted theory to be parts of data needed to specify the interpretation.
	\end{uwaga}
	
	\begin{uwaga}\label{uwaga_twarze_utb}
		For sequential theories $\USB_K$ is inter-definable with the following theory $\USB_K^{<\infty}$, which extends $U$ with all axioms of the form:
		
		\[\forall \alpha \bigl(\Ass(\phi,\alpha)\rightarrow S(\num{\qcr{\phi}}, \alpha)\equiv \phi[\alpha]\bigr),\]
		where $\Ass(\phi,\alpha)$ means that $\alpha$ is an assignment for $\phi$ and $\phi[\alpha]$ denote a formula
		\[\exists x_1,\ldots,x_n (\bigwedge_{i\leq n}\alpha(v_i) = x_i \wedge \phi[x_1/v_1,\ldots, x_n/v_n]),\]
		where $v_0,\ldots,v_n$ are all free variables of $\phi$. 
		
		Moreover, if $\Lang$ extends the language of arithmetic and $K = id$, then $\USB_K(\Lang)$ mutually syntactically definable modulo $\Lang$ with the theory $\UTB_K(\Lang)$, which extends $U$ with all axioms of the form $\forall x \bigl(T(\qcr{\phi(\dot{x})})\equiv \phi(x)\bigr),$ for $\phi(x)\in \Lang$. 
		
	\end{uwaga}

	\begin{uwaga}
		Theories $\TB$ and $\UTB$ have been studied in the literature, mostly in the context of Peano Arithmetic, $\PA$. For such a base theory, normally $\TB$ and $\UTB$ denote extensions of $\TB_{id}$ and $\UTB_{id}$ with the induction scheme extended to all formulae of $\Lang_{\TB}$ (see e.g. \cite{hal11}). Theories which are based only on the T-biconditionals are often called $\TB^-$ ($\UTB^-$) or $\TB\restr{}$ ($\UTB\restr{}$).
	\end{uwaga}

	\begin{uwaga}
		The term "definability" is used in three different contexts in the current paper. First of all, given a model $\M$, an element $a\in M$, or a subset $A\subseteq M$ can be definable in $\M$. In such contexts "definability" has its standard meaning known from basic model-theory. Secondly, a theory $U$ can be syntactically or semantically definable in another theory $V$. For the sake of precision we always speak of theories (and not models) and qualify the definability to be either semantic or syntactic. Thirdly, we consider the axiomatic theory of definability, $\DEF$.
	\end{uwaga}
	
	\section{Model-theoretical characterizations of semantical predicates}
	
	It is very easy to check the following properties of our target axiomatic theories:
	
	\begin{stwierdzenie}
		Let $U$ be any sequential $\Lang$-theory and $N: S^1_2\lhd U$.
		\begin{enumerate}
			\item $\TB_N$ uniformly imposes $\Lang$-elementary equivalence.
			\item $\DEF_N$ uniformly preserves $\Lang$-definability.
			\item $\USB_N$ uniformly imposes $\Lang$-elementarity.
		\end{enumerate}
	\end{stwierdzenie}
	
	The next theorem says that each of the above model-theoretical properties characterizes the respective theory up to syntactical definability.
	
	\begin{tw}\label{tw_impos_elem}
		Suppose that $U$ is an r.e. sequential theory in a finite language and let $N: S^1_2\lhd U$. Then
		\begin{enumerate}
			\item if $U$ uniformly imposes $\mathcal{L}$-elementarity, then $U$ syntactically defines $\USB_N(\Lang)$ modulo $\Lang$.
			\item if $U$ uniformly preserves $\Lang$-definability, then $U$ syntactically defines $\DEF_N(\Lang$) modulo $\Lang$.
			\item if $U$ uniformly imposes $\Lang$-elementary equivalence, then $U$ syntactically defines $\TB_N(\Lang$) modulo $\Lang$.
		\end{enumerate}
	\end{tw}
	
	We start with a lemma which is purely model-theoretical. We say that a formula $\phi(x)\in\Lang$ is an $\Lang$-definition if for some formula $\psi(x)$, $\phi(x):= \psi(x)\wedge \forall y\neq x \neg\psi(y)$.
	
	\begin{lemat}\label{lem_main_impos_elem}
		Let $U$ be any theory.
		\begin{enumerate}
			\item If $n\in\omega$ witnesses that $U$ uniformly  imposes $\mathcal{L}$-elementarity , then for every $\Lang$-formula $\phi(\bar{x})$ there exists a $\Sigma^*_n(\Lang_U)$ formula $\theta_{\phi}(\bar{x})$ such that $U\vdash \forall x\bigl(\theta_{\phi}(\bar{x})\equiv \phi(\bar{x})\bigr).$
			\item If $n\in\omega$ witnesses that $U$ uniformly preserves $\mathcal{L}$-definability, then for every $\Lang$-definition $\phi(x)$ there exists a $\Sigma^*_n(\Lang_U)$ formula $\theta_{\phi}(x)$ such that $U\vdash \forall x\bigl(\theta_{\phi}(x)\equiv \phi(x)\bigr).$
			\item If $n\in\omega$ witnesses that $U$ uniformly imposes $\Lang$-elementary equivalence, then for every $\Lang$-sentence $\phi$ there exists a $\Sigma^*_n(\Lang_U)$-sentence $\theta_{\phi}$ such that $U\vdash \theta_{\phi}\equiv \phi.$
		\end{enumerate} 
	\end{lemat}
	\begin{proof}
		We show the argument for $\Lang$-elementarity, the rest of them being analogous. Fix $U$, $n$ and $\phi(x)$ (without loss of generality we can assume that $\phi$ has a single free variable). If $\phi(x)$ is inconsistent with $U$, then we know what to do. Otherwise, let $\M\models U$ be such that for some $a$, $\M\models \phi(a)$. Since $n$ witnesses that $U$ imposes $\Lang$-elementarity, then
		\[U + n\mhyphen \Diag(\M)\vdash \phi(a).\]
		Consequently, there exists a $\Sigma^*_n(\Lang_U)$ sentence $\theta''_{\phi}(a,\bar{b})\in n\mhyphen \Diag(\M)$ with parameters $a,\bar{b}\in M$ such that (we assume that $a$ and $\bar{b}$ are disjoint)
		\[U\vdash \theta''_{\phi}(a,\bar{b})\rightarrow \phi(a).\]
		We may safely assume that $a,\bar{b}$ do not belong to the language of $U$. It follows that for $\theta'_{\phi}(x):= \exists \bar{y}\theta''_{\phi}(x,\bar{y})$ we have
		\[U\vdash \forall x \bigl(\theta'_{\phi}(x)\rightarrow \phi(x)\bigr).\]
	
		Let $A:= \set{\psi(x)\in\Sigma^*_n(\Lang_U)}{U\vdash \forall x(\psi(x)\rightarrow \phi(x))}$. We claim that for some $\psi_0,\ldots,\psi_k\in A$ we have
		\[U\vdash \phi(x)\rightarrow \bigvee_{i\leq k}\psi_i(x).\]
		Suppose not. Then $U+\phi(b) + \set{\neg\psi(b)}{\psi(x)\in A}$ is consistent ($b$ is a fresh constant). Let $\M$ be any model of this theory. We claim that $U+n\mhyphen \Diag(\M)+\neg\phi(b)$ is consistent as well. Indeed, for if not, then for some $\theta(x,\bar{y})$ such that $\theta(b,\bar{c})\in n\mhyphen \Diag(\M)$
		\[U\vdash \forall x\bigl(\exists \bar{y} \theta(x,\bar{y})\rightarrow \phi(x)\bigr).\]
		However, then $\exists\bar{y}\theta(x,\bar{y})\in A$, hence $\M\models \neg\exists \bar{y}\theta(b,\bar{y})$. This contradicts the fact that $\theta(b,\bar{c})\in n\mhyphen \Diag(\M)$. Hence $U+n\mhyphen \Diag(\M)+\neg\phi(b)$ is consistent, which in turn contradicts the fact that $U$ uniformly imposes $\Lang$-elementarity.
	\end{proof}
	\begin{proof}[Proof of Theorem \ref{tw_impos_elem}]
		Similarly to the case of the above lemma, we focus on $\USB_N$. Let us define the following function $\sigma:\Form_{\Lang}\rightarrow \Proof_U$:\footnote{$\Proof_U$ denotes the set of G\"odel codes of proofs from the axioms of $U$.}
		\begin{equation*}
			\sigma(\phi(x))= \min\set{p}{\exists \theta(x)\in\Sigma^*_n(\Lang_U)\Proof_U\bigl(p, \qcr{\forall x\bigl(\theta(x)\equiv \phi(x)\bigr)}\bigr)}.
		\end{equation*}
		In the above, $\Proof_U(x,y)$ is the canonical $\Sigma_1$-provability predicate for $U$. Lemma \ref{lem_main_impos_elem} makes it obvious that this function is total and recursive. For a formula $\phi(x)\in\Lang$, let $\theta_{\phi}(x)$ be the unique $\Sigma^*_n(\Lang_U)$-formula such that $\sigma(\phi(x))$ is a proof of $\forall x\bigl(\theta_{\phi}(x)\equiv \phi(x)\bigr).$ 
		The mapping $\phi\mapsto \theta_{\phi}$ is clearly total and recursive and we take $\rho(x,y)$ to be the formula which strongly represents this function in $U$. Clearly we may assume that $U$ proves that $\rho(x,y)$ is a partial function. As we have already discussed in Fact \ref{fakt_partial_truth}, there are a cut $J$ of $N$ and a satisfaction predicate $\Sat_n(s,x)$ such that $U$ proves that $\Sat_n(s,x)$ satisfies compositional clauses for all assignments (in the sense of $N$) $s$ and all $\Sigma^*_n(\Lang_U)$ formulae from $J$. Let us put:
		\begin{equation*}
			S_U(x,y):= \form_{\Lang}(x) \wedge \exists z \bigl(\rho(x,z) \wedge \Sat_n(z,y)\bigr). 
		\end{equation*}
		Within $U$ we prove that for every $\phi(v)\in\Lang$ 
		\[U\vdash \forall x \bigl(S_U(\num{\qcr{\phi(v)}},x)\equiv \phi(x)\bigr).\]
		Fix $\phi(v)$ and let $\rho(\qcr{\phi(v)})=\qcr{\theta(v)}$. Working in $U$ observe that the following equivalences hold for an arbitrary $a$
		\begin{align*}
			S_U(\qcr{\phi(v)},a) &\iff \exists z \bigl(\rho(\num{\qcr{\phi(v)}},z) \wedge \Sat_n(z,a)\bigr) \\
			&\iff \Sat_n(\num{\qcr{\theta_{\phi}(v)}}, a)\\
			&\iff \theta_{\phi}(a)\\
			&\iff \phi(a)
		\end{align*}   
		The first one follows from the definition, the fact that $\phi(v)$ is an $\Lang$-formula and the representability of the set of G\"odel codes of formulae. The second one is obtained by the representability of $\rho$. The third one holds by the appropriate version of ''"It's snowing"-It's snowing'' lemma from \cite{viss_sivs}. The last one follows, since  $\phi(v)$ is equivalent to $\theta_{\phi}(v)$ in $U$.
	\end{proof}
	
	\begin{uwaga}
		The last argument was used for the first time by Albert Visser's in the proof of \cite[Theorem 3.11]{viss19enayat}.
	\end{uwaga}
	\begin{stwierdzenie}
		There is a theory $U\supseteq \PA$ in a countable language which uniformly imposes arithmetical elementarity but does not semantically define $\TB_{id}$ modulo $\Lang_{\PA}$.
	\end{stwierdzenie}
	\begin{proof}
		We define $U$. Apart from the axioms of $\PA$, for every formula $\phi(x)$, $U$ has a binary predicate $S_{\phi(x)}$ and an axiom
		\[\forall y \bigl(S_{\phi(x)}(\num{\qcr{\phi(x)}},y)\equiv \phi(y)\bigr).\]
		To see that $U$ does not semantically define $\TB_{id}$ observe first that each model $\M\models\PA$ expands to a model of $U$ in which for every $\phi$, $S_{\phi}$ is arithmetically definable without parameters. Indeed, fixing $\M$ and  $\phi(x)\in\Lang$, one can put $S_{\phi(x)}^{\M} = \set{\pair{\qcr{\phi(x)}}{a}}{a\in M\wedge \M\models \phi(a)}$. Hence if one starts with the standard model $\mathbb{N}\models \PA$ and let $\mathbb{N}^*$ be the above defined expansion of $\mathbb{N}$ to a model of $U$, there can be no definition of $\TB_{id}$ in $\mathbb{N}^*$.
	\end{proof}
	
	Now we pass to semantical definability. The following notion will come in handy:
	
	\begin{definicja}\label{def_definable_theories}
		Let $\M$ be a sequential model for a language $\Lang$ and $N: S^1_2\lhd \ElDiag(\M)$. We say that an $\Lang_{\M}$-theory $V$ \emph{is representable in $\M$ w.r.t. $N$} (or \emph{representable in $(\M,N)$ for short})  with a formula $\theta(x,y)\in\Lang_{\M}$ iff for every $\Lang$-formula $\phi(v)$ and every $a$
		\[\phi(a)\in V \iff \M\models \theta(\num{\qcr{\phi}},a).\]
		If $V$ is representable in $\M$ with $\theta$, then we put
		\[\Proof_{\theta}(y,x,z):= \exists \psi \bigl(\theta(\psi, z)\wedge \Proof_{\emptyset}(y,\qcr{\psi\rightarrow x})).\]
		We recall that $\Proof_{\emptyset}(y,x)$ is the ($N$-translation of the) standard decidable predicate strongly representing the relation "$y$ is a proof of $x$ in pure first-order-logic". Intuitively, $\Proof_{\theta}(y,x,z)$ expresses that $y$ is a proof of $x$ which may use as an axiom any formula $\phi(x)$, such that $\phi(z)\in V$, where $V$ is the theory represented by $\theta$. 
	\end{definicja}
	
	\begin{tw}\label{tw_nonuniform_def_tb}
		If $U$ is a sequential, r.e. theory in a finite language which imposes elementary equivalence (preserves $\Lang$-definability) and $N: S^1_2\lhd U$, then $U$ parameter-free semantically defines $\TB_N$ ($\DEF_N$) modulo $\Lang$.
	\end{tw}
	\begin{proof}
		We do the case for $\DEF_N$, the proof for $\TB_N$ being fully analogous. Suppose $U$ preserves $\Lang$-definability and fix $\M$ and a suitable $n$. Then $U + n\mhyphen \Diag(\M)$ uniformly preserves $\Lang$-definability. Hence, by Lemma \ref{lem_main_impos_elem}, for every $\Lang$-definition $\phi(x)$ there is a $\Sigma^*_n(\Lang_{\M})$-formula $\theta_{\phi}(x,b)$ such that
		\[U + n\mhyphen \Diag(\M)\vdash \forall x(\theta_{\phi}(x,b)\equiv \phi(x)),\]
		So for every $\phi$ there are formulae $\psi(b)\in n\mhyphen \Diag(\M), \chi\in U$ such that
		\begin{equation*}\label{U_prov_DEF_sem}
			\emptyset\vdash  \forall y (\chi\wedge\psi(y)\rightarrow \forall x(\theta_{\phi}(x,y)\equiv \phi(x))).
		\end{equation*}
		Observe that $n\mhyphen\Diag(\M) \cup U$ is (parameter-free) representable in $(\M, N)$, thanks to the existence of partial satisfaction predicates.
		Let $\zeta$ be the representing formula. 
		For an arbitrary formula $\phi$, let $\phi^d(x):= \phi(x)\wedge \forall y\neq x \ \phi(y)$. We write the definition $D_{\M}(\phi, a)$ in $\M$:
			\begin{multline*}
				\exists \theta(x,y)\exists p\exists b\bigl[\Sat_n^*(\theta,\pair{a}{b})\wedge  \forall y \neq a \neg \Sat_{n}^*(\theta,\pair{y}{b})\wedge  \Proof_{\zeta}(p, \qcr{\forall x (\theta(x,y)\equiv \phi^d(x))}, b) \wedge \\ \wedge \forall p'<p\forall \theta'<p\bigl(\theta'\in\Sigma_n^*(\Lang_U)\rightarrow \forall b\neg\Proof_{\zeta}(p', \qcr{\forall x (\theta'(x,y)\equiv \phi^d(x))},b)\bigr) \bigr].
			\end{multline*}	
		We prove that $D_{\M}(x,y)$ defines $\DEF_N$ in $\M$. Fix $\phi(x)$ and fix the least proof $p$ with the following property
			\begin{center}
				There are $b\in M$, $\chi\in U$ and $\psi(y),\theta(x,y)\in \Sigma^*_n(\Lang_U)$ such that $M\models \psi(b)$ and it holds that 
				$\Proof_{\emptyset}(p, \qcr{\forall y (\chi\wedge\psi(y)\rightarrow \forall x(\theta(x,y)\equiv \phi^d(x)))}).$
			\end{center}
		$p$ is well-defined by the well-foundedness of natural numbers. Let $\qcr{p}$ denote the G\"odel code of $p$. By definition
		\begin{multline*}
			\M\models \exists b\exists \theta \in\Sigma_n^*(\Lang_U)\Proof_{\zeta}(\num{\qcr{p}}, \qcr{\forall x (\theta(x,y)\equiv \phi^d(x))}, b) \wedge \\ \wedge \forall p'<\num{\qcr{p}}\forall \theta'<\num{\qcr{p}}\bigl(\theta'\in\Sigma_n^*(\Lang_U)\rightarrow \forall b\neg\Proof_{\zeta}(p', \qcr{\forall x (\theta'(x,y)\equiv \phi^d(x))},b)\bigr).
		\end{multline*}
		Moreover $\num{\qcr{p}}\in M$ is the unique, up to the $E$-equivalence relation, element of the domain of $\M$ for which the above holds. Assume now $\M\models D_{\M}(\num{\qcr{\phi}},a)$. Then clearly for some $b$ and standard $\theta\in\Sigma_n^*(\Lang_U)$, $\M\models \Sat_n^*(\theta,\pair{a}{b})\wedge  \forall y \neq a \neg \Sat_{n}^*(\theta,\pair{y}{b})$. Moreover, for the same $b$, $\theta$, $U+n\mhyphen\Diag(\M)\vdash  \forall x (\theta(x,b)\equiv \phi^d(x))$, hence in $\M$, $\phi(x)$ defines $a$. Now assume that $\phi(x)$ defines $a$ in $\M$. Then obviously $\phi^d(x)$ defines $a$. Fix any $b$ and $\theta$ such that $\M\models \Proof_{\zeta}(\qcr{p}, \qcr{\forall x (\theta(x,y)\equiv \phi^d(x))}, b)$. By definition and the fact that both $p$ and $\theta(x,y)$ are standard,
		it follows that $\M\models \forall x (\theta(x,b)\equiv \phi^d(x))$. Hence $\theta(x,b)$ defines $a$ in $\M$ and $D_{\M}(\num{\qcr{\phi}},a)$ follows by ""It's snowing"-It's snowing" lemma for $\Sat_n$. 
	\end{proof}

	The lemma below provides a criterion for definability of $\USB_N$ in a sequential model. It will be used to characterize the satisfaction predicate up to semantical definability and, in the next section, separate satisfaction from truth and definability.

	\begin{lemat}\label{lem_main_model_theor_def}
		Let $\M$ be a model of a sequential theory
		for some finite language $\Lang'$ and $N: S^1_2\lhd \ElDiag(\M)$. Then $\USB_N(\Lang)$ is syntactically definable in $\ElDiag(\M)$ modulo $\Lang$ iff for some $(\M,N)$-representable, sequential $\Lang'_{\M}$-theory $V$ consistent with $\ElDiag(\M)$, $V\vdash \ElDiag(\M\restr{\Lang})$.
	\end{lemat}
	\begin{proof}
		Assume that $\ElDiag(\M)$ defines $\USB_N(\Lang)$ via a formula $S'(\bar{x})$ and let $n$ be the complexity of $S'$. Let $V$ be the theory $\USB_N(\Lang)[S'/S] + n\mhyphen \Diag(\M)$. This theory is representable in $(\M,N)$, by the existence of partial satisfaction predicates. Then clearly, for every $\phi(x)\in\Lang$ and all $a\in M$ we have 
		\begin{align*}
			\M\models\phi(a)&\iff \M\models S'(\num{\qcr{\phi(x)}}, a)\\
			&\iff S'(\num{\qcr{\phi(x)}}, a)\in n\mhyphen\Diag(\M)\\
			&\iff V\vdash \phi(a).
		\end{align*}

		Assume now, $V\vdash \ElDiag(\M\restr{\Lang})$. It follows from this, the consistency of $V$ and the completeness of $\ElDiag(\M)$ that for every $\phi(a)\in\Lang_{\M}$, $\phi(a)\in\ElDiag(\M)$ iff there is $\psi(x)\in \Lang'$ such that $\psi(a)\in V$ and
		\[\emptyset\vdash \psi(x)\rightarrow \phi(x).\]
		Let $\theta$ represent $V$ in $\M$. Now, we can define $S_{\M}(x, c)$ as
		\[\exists y\bigl(\Proof_{\theta}(y,x,c) \wedge \forall z \bigl(z< y \rightarrow \neg \Proof_{\theta}(z,x,c)\bigr)\bigr) \]
		We check that this definition works, following the lines of the standard Rosser-style argument. Fix $\phi(v)\in\Lang$ and any  $a\in M$. First assume that $\phi(a)$ holds in $\M$. As above for some $\Lang'$-formula $\psi(v)$ such that $\psi(a)\in V$ 
		\[\emptyset\vdash \psi(v)\rightarrow \phi(v).\]
		Let $l\in \omega$ code this proof. Clearly then $\M\models \Proof_{\theta}(\num{l}, \num{\qcr{\phi(v)}}, a)$. We claim that for no $l'<l$ and for no $\xi$ such that $\M\models\Sat_{n}^*(\xi,a)$ we have $\M\models\Proof_{\theta}(\num{l'},\qcr{\xi\rightarrow\neg\phi(v)}, a)$.	
		Indeed, if this were the case, then $\xi$ would be a standard formula. Since $\M\models V$, we would have $\M\models \xi(a)$. However, obviously also $\xi(a)\rightarrow\neg\phi(a)$ is true in $\M$. This contradicts $\phi(a)$ being true in $\M$. Assume now that $\neg\phi(a)$ holds in $\M$. Then, arguing as previously we have that for some $l\in\omega$, $\M\models \Proof_{\theta}(\num{l}, \num{\qcr{\neg\phi(v)}}, a)$ and for all $l'<l$, $\M\models\neg \Proof_{\theta}(\num{l'}, \num{\qcr{\phi(v)}}, a)$. Hence, 
		\[\M\models \forall y\bigl(\Proof_{\theta}(y,(\num{\qcr{\phi(v)}},c) \rightarrow \exists z \bigl(z< y \wedge \Proof_{\theta}(z,(\num{\qcr{\neg\phi(v)}},c)\bigr).\]
		Hence $\M\models \neg S_{\M}(\num{\qcr{\phi}},a)$.
	\end{proof}
	
	The theorem below is a semantical variant of Theorem \ref{tw_impos_elem} (for $\USB$). We stress that neither of the theorems directly implies the other one. 
	\begin{tw}\label{tw_nonuniform_impos_elem_utb}
			If $U$ is an r.e. sequential theory in a finite language $\Lang'$ which imposes $\Lang$-elementarity, then $U$ semantically defines $\USB_N(\Lang)$ modulo $\Lang$ for each interpretation $N:S^1_2\lhd U$.
	\end{tw}
	\begin{proof}
		Fix $U$, $N$ and let $\M\models U$. Since $U$ imposes $\Lang$-elementarity, there is an $n$ such that $U+n\mhyphen\Diag(\M)\vdash\ElDiag(\M\restr{\Lang})$. Since $U$ is r.e., sequential and in a finite language, $U+ n\mhyphen\Diag(\M)$ is representable in $(\M,N)$. Hence by Lemma \ref{lem_main_model_theor_def}, $\M$ defines $\USB_N({\Lang})$.
	\end{proof}
	
	\section{Definability between axiomatic theories of truth, definability and satisfaction}
	
	\label{sec_definability_between_theories}
	
	Now, we turn to the definability relations between the theories of semantical notions in study. We say that $U$ has a definable element, if for some formula $\phi(x)$, $U\vdash \exists ! x \phi(x)$.
	
	\begin{stwierdzenie}\label{stw_latwa_definiowalnosc}
		For every sequential theory $U$ and interpretation $N: S^1_2\lhd U$, $\TB_N \leq_{\Lang_U} \USB_N$ and  $\DEF_N \leq_{\Lang_U}\USB_N$. Moreover, if $U$ has a definable element, then $\TB_N\leq_{\Lang_U}\DEF_N$.
	\end{stwierdzenie}
	\begin{proof}
		The definability $\TB_N \leq_{\Lang_U} \USB_N$ is trivial. Moreover, clearly
		\[D(\phi(x), y):= \exists! z S(\phi(x), z)\wedge S(\phi(x),y)\]
		witnesses the definability of $\DEF_N$ in $\USB_N$. For the remaining case, let $\theta(x)$ be a $U$-provable definition. Then, the formula
		\[T(\phi):= \exists y D(\qcr{\theta(x) \wedge \phi}, y),\]
		witnesses the definability of $\TB_N$ in $\DEF_N$. In the above, $\qcr{\theta(x)\wedge \phi}$ denotes the ($S^1_2$-representation of the) function with one free variable $\phi$ which, given a (G\"odel code of a) sentence $\phi$ returns the (G\"odel code of the) sentence $\theta(x)\wedge \phi$. This function need not be provably total in $U$, however for every standard sentence $\phi$ the value of the function exists and this all we need for the definability of $\TB_N$.
	\end{proof}
	
	In general, it is not true that $\DEF_N$ syntactically defines $\TB_N$ modulo $U$. This is due to the following simple observation:
	
	\begin{stwierdzenie}\label{stw_def_no_tb}
		If $U$ has a model in which no element is parameter-free definable, then $\DEF_N$ does not semantically define $\TB_N$ modulo $\Lang_U$.
	\end{stwierdzenie}
	\begin{proof}
		If $\M\models U$ and has no parameter-free definable elements, then $\DEF_N$ is definable in $\M$ via a parameter-free formula $\phi(x,y):= x\neq x \wedge y\neq y$. Clearly, $\TB_N$ cannot be parameter-free definable in $\M$, by Tarski undefinability theorem for $\Th(\M)$.
	\end{proof}
	
	A concrete example of a sequential theory with a model with no parameter-free definable elements is AS. In particular,  $(V_{\omega}, \in)\sqcup (V_{\omega}, \in)$  is a model of AS in which no element is definable (for models $\M$ and $\N$, $\M\sqcup\N$ denotes their disjoint union), because every element can be moved by an automorphism.

	We proceed to separations between $\TB_N$, $\USB_N$ and $\DEF_N$. In the specific context where $U$ is a theory in an arithmetical language consistent with $\PA$ and $N$ is the identity interpretation, the separation between $\TB_N$ and $\USB_N$ was first shown in \cite{wcislyk_models_weak}. Below we prove that
	$\USB_N$ is not semantically definable in $\DEF_N$ which obviously generalizes also to the case of $\TB_N$.

	\begin{tw}\label{tw_nondef_sat_in_def_tb}
		Suppose $U$ is a consistent, sequential r.e. $\Lang$-theory and $N: S^1_2\lhd U$. Then neither $\DEF_N(\Lang)$ nor $\TB_N(\Lang)$ semantically defines $\USB_N(\Lang)$ modulo $\Lang$.
	\end{tw}
	\begin{proof}
		We do the case of $\DEF_N(\Lang)$, the argument for $\TB_N(\Lang)$ being analogous. Fix $U,N$. Let $\M\models U$ be any $\Lang$-recursively saturated model.

		By recursive saturation there is $c\in M$ such that for any $\phi(x)$ and any $a\in M$
		\[\M\models \pair{\num{\qcr{\phi(x)}}}{a}\in_{\M} c \equiv \exists!x\phi(x)\wedge \phi(a).\]
		In the above $\pair{\cdot}{\cdot}$ denotes the pairing function given by the direct interpretation of AS. Let $D^{\M}:= \set{\pair{a}{b}}{\M\models \pair{a}{b}\in_{\M} c}$. Then $(\M,D^{\M})\models \DEF_N(\Lang)$ and $D^{\M}$ is definable with a parameter in $\M$. Hence there cannot be any definition of $\USB_N(\Lang)$ in $(\M, D^{\M})$ by Tarski's undefinability of truth theorem (applied to $\ElDiag(\M)$). 
	\end{proof}
	
	Now, we shall separate $\DEF_N$ from $\TB_N$. Unlike in the case of $\DEF_N$ and $\USB_N$,
	 we cannot count for the undefinability result over all sequential theories, since some of them trivialize the behaviour of $\DEF_N$.
	\begin{przyklad}
		Define $\textsc{AS}_f$ to be the extension of AS with a sentence "$f$ is an automorphism which moves every element". Then, by the remark following Proposition \ref{stw_def_no_tb}, $\textsc{AS}_f$ is consistent. Clearly in every model $\M\models \textsc{AS}_f$, $f^{\M}$ is an automorphism of $\M$ which moves every element. Hence $\DEF_N$ is syntactically definable in $\textsc{AS}_f$ via any formula $\phi(x,y)$ with provably empty extension.
	\end{przyklad} 
	

As we already noted in Proposition \ref{stw_latwa_definiowalnosc},
in the case when $U$ admits a provably definable element, $\TB_N$ is contained in $\DEF_N$. We will shortly see, in Theorem \ref{tw_nondef_def_in_tb}, that under relatively mild assumptions this containment is in fact strict.  However, if $U$ has provably finitely many definable elements, then $\DEF_N$ does not exceed $\TB_N$:
	
	\begin{stwierdzenie}
		If for finitely many formulae $\theta_0,\ldots, \theta_n$, $U$ proves $\bigwedge_{i\leq n}\exists ! \theta_i(x)$ and for an arbitrary formula $\phi(x)$
		\[U\vdash \exists ! \phi(x)\rightarrow \forall x \bigl(\phi(x)\rightarrow \bigvee_{i\leq n}\theta_i(x)\bigr),\]
		then $\TB_N$ syntactically defines $\DEF_N$.
	\end{stwierdzenie}
	\begin{proof}
		Fix $\was{\theta_i}$ as above and working in $\TB_N$, put
		\[D(\psi(x),y):= \Form_U(\psi(x)) \wedge \bigvee_{i\leq n} \bigl[T(\qcr{\forall x (\psi(x)\equiv \theta_i(x))})\wedge \theta_i(y)\bigr].\]
		In the above, as in the proof of Proposition \ref{stw_latwa_definiowalnosc},  $\qcr{\forall x (\psi(x)\equiv \theta_i(x))}$ denotes the ($S^1_2$-representation of the) function which, given a (G\"odel code of a) formula $\psi(x)$ with a single free variable returns the (G\"odel code of the) sentence $\forall x (\psi(x)\equiv \theta_i(x))$.
	\end{proof}
	
	Also, it is fairly easy to convince oneself that for each $n$ there is a sequential theory which admits exactly $n$ definable elements. We shall show that if $U$ is a sequential, r.e. theory in a finite language $\Lang$ which for some $n$ has a model with infinitely many $\Sigma^*_n(\Lang)$-definable elements then $\TB_N$ does not semantically define $\DEF_N$ modulo $\Lang$. Let us observe that each theory $U$ that admits a consistent extension which directly interprets ASE via a one-dimensional and parameter-free interpretation is of this form. 
	
	The idea of the proof is to construct a model $\M\models U$  with $c\in M$ such that 
	
	\begin{enumerate}
		\item $c$ codes the theory of $\M$,
		\item there is a partition $\was{M_i}_{i\in \omega}$ of $\M$ such that for every $n$ there exists an element $d_n\in \Def(\M)$ which is not $\Sigma_n^*(\Lang_{M_n})$-algebraic in $\M$.
	\end{enumerate}
	This will ensure that for any $n$, the $\Lang_T$ model $(\M,c^{\M})$ has an $\Lang_T$-elementarily equivalent, $\Sigma_n^*(\Lang_T)$-elementary extension which differs from $\M$ on definable elements. Moreover, $\ElDiag(\M)$ does not syntactically define $\DEF_N$. In our construction we shall make use of flexible formulae.
	
	\begin{definicja}
		Let $U$ be any theory in a language $\Lang$. A formula $\theta(x)$ is $(n,k)$-flexible for $U$ iff for every $\M\models U$ and every $\Sigma^*_k(\Lang)$ formula $\delta(x)$ there exists $\M\preceq_{n}\N$ such that $\N\models \theta(x)\equiv \delta(x)$. 
	\end{definicja}
	
	The following is a variation of \cite[Theorem 11, Section 2.3]{lindstrom_aspects}.
	
	\begin{stwierdzenie}\label{prop_flexible}
		If $U$ is a sequential, r.e. theory in a finite language, then for every $n,k$ there exists an $(n,k)$-flexible formula for $U$.
	\end{stwierdzenie}
	\begin{proof}
		Fix $U$ and $n,k\in\omega$. Let $\Proof_{\theta}(y,x,z)$ be the provability predicate defined as in Definition \ref{def_definable_theories} for $\theta(x,z):= \Sat_{n+1}(x,z)\vee U(x)$ (where $U(x)$ is any decidable predicate which represents an axiom set of $U$ and $\Sat_{n+1}$ is a partial truth predicate for $\Lang_U$). Let $\rho'(y, \phi,\psi)$ be the formula
		\[\psi \in\Sigma^*_k(\Lang_U)\wedge\exists z\Proof_{\theta}(y, \qcr{\neg \forall x \bigl(\phi\equiv \psi\bigr)}, z),\]
		In the above, $\qcr{\neg \forall x \bigl(\phi\equiv \psi\bigr)}$ denotes the ($S^1_2$-representation of the) function which, given  (G\"odel codes of) formulae $\phi(x),\psi(x)$ (both having at most $x$-free) returns the (G\"odel code of the) sentence $\neg\forall x (\phi\equiv \psi)$. Finally put 
		\[\rho(y, \phi,\psi):=\rho'(y, \phi,\psi )\wedge \forall y'\forall \psi'\bigl(\pair{y'}{\psi'}<\pair{y}{\psi}\rightarrow \neg\rho'(y', \phi,\psi')\bigr).\]
		Using the diagonal lemma in $U$ define $\gamma(x)$ to be the formula such that
		\[U\vdash \gamma(x)\equiv \exists y\exists\psi\bigl(\rho(y,\num{\qcr{\gamma}},\psi)\wedge \Sat_{k}(\psi,x)\bigr).\]
		Now fix any $\M\models U$, any $\psi\in\Sigma^*_k(\Lang_U)$ and aiming at a contradiction assume that there is no $\Sigma^*_n(\Lang_U)$-elementary extension of $\M$ making $\forall x \bigl(\gamma(x)\equiv \psi(x)\bigr)$ true. It follows that $U+ (n+1)\mhyphen\Diag(\M)\vdash \neg \forall x \bigl(\gamma(x)\equiv \psi(x)\bigr).$ Hence, for some $\chi\in (n+1)\mhyphen\ElDiag(\M)$
		\[U\vdash \chi\rightarrow \neg \forall x \bigl(\gamma(x)\equiv \psi(x)\bigr).\]
		By existentially generalizing on the additional parameters from $\M$ we can assume that $\chi$ is in the language of $U$. Let us pick $\psi\in M$ and the proof $p\in M$ such that $\pair{p}{\psi}$ is minimal such that $\M\models\rho'(\num{\qcr{p}}, \num{\qcr{\gamma}},\num{\qcr{\psi}})$. Such a pair exists, by the least number principle in the well-founded part of $\M$. In particular, $p$ witnesses that for some $\chi'\in (n+1)\mhyphen \Diag(\M)$, $U\vdash \chi'\rightarrow \neg\forall x (\gamma \equiv \psi)$. As previously we can assume that $\chi'$ contains no parameters. Then obviously $\M\models \rho(\num{\qcr{p}},\num{\qcr{\gamma}}, \num{\qcr{\psi}})$. Moreover, for any $y$ and $z$, if $\M\models \rho(y,\num{\qcr{\gamma}}, z)$, then $y$ is in the $E$-equivalence class corresponding to $\qcr{p}$
		and $z$ is in the $E$-equivalence class corresponding to $\qcr{\psi}$ (See Convention \ref{konwencja_arytmetyka_w sekwencji}). 
		
		It follows that $\M\models \forall x \bigl(\gamma(x)\equiv \Sat_{k}^*(\num{\qcr{\psi}},x)\bigr)$. Consequently, $\M\models \forall x \bigl(\gamma(x)\equiv \psi(x)\bigr)$. However, since $\M\models U + \chi'$, $\M\models \neg\forall x\bigl(\gamma(x)\equiv \psi(x)\bigr)$. 
	\end{proof}
	
	\begin{lemat}\label{lem_hierarchical_model}
		For every r.e. sequential theory $U$ in a finite language $\mathcal{L}$, every $n,k\in\omega$ and every model $\M\models U$ which contains infinitely many $\Sigma^*_k(\Lang)$-definable elements there is an $\N\succeq_n\M$ and a decomposition $N = \bigcup_{m\in\omega} M_m$ such that for every $m,l\in\omega$ there is an element of $N$ which is parameter-free definable in $\N$ but not $\Sigma^*_l(\Lang_{M_m})$-algebraic.
	\end{lemat}
	\begin{proof}
		Fix $U$, $\M$ and $n,k$ as above. Without loss of generality assume that $k<n$. We shall define 
		\begin{itemize}
			\item a chain of models $\was{\M_i}_{i\in\omega}$;
			\item a sequence of formulae $\was{\theta_i(x)}_{i\in\omega}$;
			\item a sequence of numbers $\was{n_i}_{i\in\omega}$. 
		\end{itemize} 
		such that
		\begin{enumerate}
			\item $n_0 = n+1$ and for each $i$, $n_i<\Sigma^*(\theta_{i+1})<n_{i+1}$, where $\Sigma^*(\theta)$ denotes the $\Sigma^*_n$ complexity of $\theta$.
			\item $\M_0 = \M$ and for each $i$, $\M_i\preceq_{n_i}\M_{i+1}$ and $\M_i\models U$.
			\item for each $i$, $\M_{i+1}$ contains an element which is definable in $\M_{i+1}$ with $\theta_{i+1}(x)$ but not $\Sigma^*_{n_i}(\Lang_{M_i})$-algebraic (in $\M_{i+1}$).
		\end{enumerate}
		Assume that we have succeeded. Let $\N = \bigcup_{i\in\omega}\M_i$. Then it follows that for every $i\in\omega$, $\M_i\preceq_{n_i}\N$. Hence, the following follows
		\begin{itemize}
			\item $\N\models U$ and $\M\preceq_n \N$.
			\item if $d$ is definable in $\M_{i+1}$ by the formula $\theta_{i+1}$, then $d$ is definable with $\theta_{i+1}$ in $\N$ and, in $\N$, $d$ is not $\Sigma^*_{n_i}(\Lang_{\M_{i}})$-algebraic.
		\end{itemize}
		Now we show how to construct the three sequences. For $n=0$ we do what is necessary
		with $n_0$ and $\M_0$ and let $\theta_0(x)$ be any formula. Assume that our sequences are constructed up to $k$. We let $\theta_{k+1}$ be an $(n_k+1, n_k+1)$-flexible
		formula. Let $\Th:= A\cup B \cup C$ be a theory in the language $\Lang(\M_{k})\cup\was{d}$, where $d$ is a fresh constant and $A,B,C$ are defined:
		\begin{itemize}
			\item[$A.$]  $(n_k+1)\mhyphen\Diag(\M_k)$
			\item[$B.$] $\set{\exists^{\leq i}x\phi(x)\rightarrow \neg \phi(d)}{\phi(x)\in \Sigma_{n_k}^*(\Lang_{\M_k}), i\in \omega}$
			\item[$C.$] $\exists ! x \theta_{k+1}(x) \wedge \theta_{k+1}(d).$
		\end{itemize}
		We claim that $\Th$ is consistent. Indeed, fix any finite fragment $K$ of $Th$,
		let $\phi_0,\ldots,\phi_m$ be all $\Sigma^*_{n_k}(\Lang_{\M_{k}})$-formulae that occur in $B\cap K$ and without loss of generality assume that each of them is satisfied by finitely many elements in $\M_k$. Pick any $\psi(x)\in\Sigma^*_k(\Lang_U)$ which defines in $\M_k$ an element which is different from any element satisfying $\bigvee_{j\leq m}\phi_j(x)$. I.e. $\psi(x)$ is such that
		\[\M_k\models  \exists !x\psi(x) \wedge \forall x \bigl(\psi(x)\rightarrow \bigwedge_{j\leq m} \neg\phi_j(x)\bigr).\] 
		Such a $\psi$ exists because in $\M$ there are infinitely many $\Sigma^*_k(\Lang_U)$ definable elements and $\M\preceq_{k} \M_k$. By flexibility of $\theta_{k+1}$ there exists $\N\succeq_{{n_k}}\M$ such that $\N\models U+ \forall x \bigl(\theta_{k+1}(x)\equiv \psi(x)\bigr)$. By the choice of $\psi$ and $n_k$, $\N\models \forall x \bigl(\psi(x)\rightarrow\bigwedge_{j\leq m}\neg\phi_j(x)\bigr)$. Hence $\N\models K$ and $\Th$ is indeed consistent.
		
		We put $\M_{k+1}$ to be any model of $\Th$. We set $n_{k+1}:= \Sigma^*(\theta_{k+1})+1$.
	\end{proof}

	\begin{tw}\label{tw_nondef_def_in_tb}
		Suppose that $U$ is a sequential, r.e. theory in a finite language $\Lang$, $K: S^1_2\lhd U$
		and $n$ is such that there is a model $\M\models U$ with infinitely many $\Sigma^*_n(\Lang)$-definable elements. Then $\TB_K(\Lang)$ does not semantically define $\DEF_K(\Lang)$ modulo $\Lang$ and does not impose the equality of definables.
	\end{tw}
	\begin{proof}
		Fix $U$, $K$, $n\in\omega$ and let $\M\models U$ have infinitely many $\Sigma^*_n(\Lang)$-definable elements. Without loss of generality assume that it is $\omega$-saturated, that is, $\M$ realizes all the types with finitely many parameters from $\M$. Then, for every $V$, complete and consistent extension of $U$, there is $c\in M$ which codes $V$ in $\M$ (according to the Definition \ref{defi_kodowanie}). 
		Let $\N$ be the $\Sigma^*_k(\Lang)$-elementary extension of $\M$ constructed as in Lemma \ref{lem_hierarchical_model}, where $m$ is greater than the complexity of $x\in_U y$ (see Definition \ref{defi_sekwencyjność}). Since $m$ is large enough, in $\N$ there is an element which codes $\Th(\N)$. 
		Hence we can fix $c\in N$ such that 
		\[(\N,c)\models \TB_{K}(\Lang).\] 
		\paragraph{Claim} For all sufficiently large $k$, for all $b\in N$ there is an $\N'\succeq_k\N$ such that 
		\begin{enumerate}
			\item $(\N',c,b)\models \Th((\N,c,b))$
			\item there is a definable element in $\N$ which is not definable in $\N'$.
		\end{enumerate}
		Before proving the claim we show that it implies that $\TB_K$ neither imposes the equality of ($\Lang$-)definables
		 nor semantically defines $\DEF_K$. Indeed, assuming the claim is true, for every $b\in \N$, $\Th((\N,c,b))$ is a complete and consistent theory which syntactically defines $\TB_K$ and does not impose the equality of definables. It directly follows that $\TB_K$ does not impose the equality of definables. Finally, if there were a definition of $\DEF_K$ in $(\N,c)$,
		 then for some $b\in N$ this definition would work provably in $\Th((\N,c,b))$, hence $\DEF_K$ would be syntactically definable in $\Th((\N,c,b))$. 
		  This would contradict the fact that $\Th((\N,c,b))$ does not impose the equality of definables .

		Fix any $k$ which is larger than the complexity of $x\in y$ and any $b\in N$. Now, aiming at a proof of the claim, we define a model $\N'$ satisfying $1$ and $2$.  Let $l$ be such that $b\in \M_l$,
		where $\was{\M_i}_{i\in\omega}$ is the chain constructed for $N$ in Lemma \ref{lem_hierarchical_model}. Now, let $d$ be any element of $N$ which is parameter-free definable, but not $\Sigma_{k+1}^*(\Lang_{M_l})$-algebraic.
		
		Consider the theory $\Th$ which is the union of the following sets (we say that a formula $\phi(x)$ is syntactically a definition if for some $\phi'(x)$, $\phi(x)$ is of the form $\phi'(x) \wedge \forall y\bigl(\phi'(y)\rightarrow y=x\bigr)$):
		\begin{itemize}
			\item $\Th((\N,c,b))$
			\item $(k+1)\mhyphen \Diag(\N)$
			\item $\set{\neg\psi(d)}{\psi(x)\in\Lang \wedge \psi \textnormal{ is syntactically a definition}}$
		\end{itemize} 
		Clearly any model of $\Th$ will satisfy 1 and 2, so it is enough to show that $\Th$ is consistent. Suppose not. Then for finitely many formulae $\psi_0,\ldots,\psi_t$
		 and some formula $\nu(x)\in\Sigma^*_{k+1}(\Lang\cup\was{c,b})$ such that $\nu(d)\in (k+1)\mhyphen\Diag(\N)$ and we have
			\[\Th((\N,c,b))\vdash \nu(x)\rightarrow \bigvee_{j\leq t}\psi_j(x).\]
			It follows that $\Th((\N,c,b))\vdash \exists^{\leq t}x\nu(x)$. However, this is impossible since $\N\models \Th((\N,c,b))$ and $d$ is not $\Sigma^*_{k+1}(\Lang\cup\was{c, b})$-algebraic in $\N$. 
	\end{proof}
	
	\begin{wniosek}\label{cor_nonres_def}
		Let $U$ be a sequential r.e. theory in a finite language with $K: S^1_2\lhd U$. Suppose that for some $n\in\omega$ there is a model $\M\models U$ with infinitely many $\Sigma^*_n(\Lang)$-definable elements, then $\DEF_K$ does not have a restricted axiomatization modulo $U$.  
	\end{wniosek}
	\begin{proof}
		Fix $U$ as above and assume $V$ is a $\Sigma_k^*$-restricted $\Lang \cup \was{D}$-theory such that 
		\[U+V \equiv \DEF_K(\Lang).\]
		Take any model $\M\models U$ and $n\in\omega$ such that in $\M$ there are infinitely many $\Sigma^*_n(\Lang)$-definable elements. Without loss of generality assume that $\M$ is recursively saturated. Hence, there is an element $c\in M$ such that for all $a\in M$, $\phi(x)\in\Lang$,
		\[\M\models \pair{\phi(x)}{a}\in_{\M} c \equiv \exists !x \phi(x)\wedge \phi(a).\]
		In the above both $\pair{\cdot}{\cdot}$, $\in_{\M}$ denote the set theoretical operations given by the direct interpretation of AS.
		 It follows that $\ElDiag(\M)\vdash \DEF_K[\pair{x}{y}\in c / D(x,y)]$. Suppose that $l$ is the complexity of $\pair{x}{y}\in_{\M} z$. By our assumption, it follows that $U+(k+l)\mhyphen \Th((\M,c))\vdash V[\pair{x}{y}\in_{\M} c / D(x,y)]$. Hence 
		\[U+(k+l)\mhyphen \Th((\M,c))\vdash \DEF_K[\pair{x}{y}\in_{\M} c / D(x,y)].\]
		Let $\N\succeq_{k+l} \M$ be the model constructed in Lemma \ref{lem_hierarchical_model}. It follows that $\N\models U+ (k+l)\mhyphen \Th(\M,c)$, hence, by our assumption
		\[\N\models \DEF_K[\pair{x}{y}\in_{\M} c / D(x,y)].\] Hence for some $m\in\omega$ in $\N$ all $\Lang$-definable elements are $\Sigma^*_m(\Lang)$-definable with a parameter $c$.
		This contradicts the properties of $\N$. 	
	\end{proof}

	We finish this section with some observations regarding the relations between semantical and syntactical definability of $\TB$, $\DEF$ and $\USB$.
	
	\begin{stwierdzenie}\label{prop_UTB_pf}
		For any sequential model $\M$ for a finite language $\Lang$ and an interpretation $N:S^1_2\lhd \ElDiag(\M)$, if $\ElDiag(\M)$ syntactically defines $\USB_N(\Lang)$, then $\Th(\M)$ syntactically defines $\USB_N(\Lang)$.
	\end{stwierdzenie}
	\begin{proof}
		Fix $\M, N$ and suppose $T'_{\M}(x,y,z)$ is such that for some $b$, $T'_{\M}(x,y,b)$ defines $\USB_N(\Lang)$ in $\M$. Let $n$ be the $\Delta^*$ complexity of $T'_{\M}$. Let $\Sat_n(x,y)$ be the partial satisfaction predicate for $\Delta_n$ formulae (as in Fact \ref{fakt_partial_truth} and the discussion following it). Let $\Sat_n'(\psi,z,w_0,w_1)$ be the formula
        \[\Form^{v_0,v_1,v_2}_{\Lang}(\psi)\wedge\Sat_n(\psi,\left[v_0\mapsto w_0,v_1\mapsto w_1,v_2\mapsto z\right]).\]
  In the above, $\Form^{v_0,v_1,v_2}_{\Lang}(\psi)$ expresses that $\psi$ is a formula whose free variables are exactly $v_0,v_1,v_2$ ($v_0,v_1,v_2$ being fixed variables) and $[v_0\mapsto x,v_1\mapsto y,v_2\mapsto z]$ denotes the unique assignment sending $v_0$ to $x$, $v_1$ to $y$ and $v_2$ to $z$. Now put
  \[\Theta(\phi,\psi,x,z):= \SAT^{w_0,w_1}_n[\Sat'_{n}(\psi, z,w_0,w_1)]\rightarrow \Sat_{n}(\psi,[v_0\mapsto \phi, v_1\mapsto x, v_2\mapsto z]).\]
	Recall that $\SAT_n^{w_0,w_1}[\xi]$ is a formula expressing that $\xi(w_0,w_1)$ is a $\Delta_n^*$-satisfaction predicate with $w_0$ being the variable corresponding to formulae and $w_1$ being the variable corresponding to assignments (as introduced in the remarks following Fact \ref{fakt_partial_truth}). Therefore, the intuitive meaning of $\Theta(\phi,\psi,x,z)$ is: if $\psi(z,w_0,w_1)$ is a satisfaction predicate for $\Delta_n^*$ formulae (where $z$ is treated as a parameter), then $\phi$ is satisfied by $x$ according to this satisfaction predicate."  Finally put $T_{\M}(\phi,x):= \forall \psi(v_0,v_1,v_2)\in\Delta^*_n \forall z\bigl(J_n(\pdpt(\psi))\rightarrow \Theta(\phi,\psi,x,z)\bigr).$

	
		The meaning of the definition is that $\phi$ is satisfied by $x$ if and only if $\phi$ is satisfied by $x$ according to all $\Sigma^*_{n}(\Lang)$ truth definitions for formulae of complexity at most the complexity of $\phi$. $T_{\M}$ clearly does not involve any parameters. Such a definition does the job since first of all there is $\Delta_n^*$ definable satisfaction predicate (i.e. $T'_{\M}(x,y,b)$) and secondly, all satisfaction predicates have to agree on formulae of standard complexity. We omit the elementary details.
	
	\end{proof}

	\begin{wniosek}\label{wniosek_ogr_def_utb}
		Let $U$ be an r.e. sequential theory in a finite language and $N: S^1_2\lhd U$. If there is an $n$ such that for every $\M\models U$, $\ElDiag(\M)$ syntactically defines $\USB_N$ with a formula of complexity $n$, then $U$ syntactically defines $\USB_N$.
	\end{wniosek}
	\begin{proof}
		It is enough to observe that the formula $T_{\M}$ from the proof of Proposition \ref{prop_UTB_pf} depends only on $n$, but not on $\M$.
	\end{proof}
	
	\begin{uwaga}
		In Theorem \ref{th_theory_with_rec_sat_without_utb} we shall construct a theory $U\supseteq\PA$ in an infinite language such that for every $\M\models U$ $\ElDiag(\M)$ syntactically defines $\UTB_{id}(\Lang_{\PA})$ with an atomic formula but $U$ does not syntactically define $\UTB_{id}(\Lang_{\PA})$. Similarly, in Theorem \ref{th_theory_finlang_with_rec_sat_without_utb}, we will construct a theory $V$ in a finite language such that in every model $\M \models V$, $\ElDiag(\M)$ syntactically defines $\UTB_{id}(\Lang_{\PA})$ without a uniform bound on the complexity of the definitions, such that $V$ does not syntactically define
	 $\UTB_{id}(\Lang_{\PA})$.  
	\end{uwaga}
	Corollary \ref{wniosek_ogr_def_utb} fails dramatically in the context of $\TB_N$.
	
	\begin{stwierdzenie}\label{prop_ogr_sem_tb_2}
		If $U$ is any sequential theory and $N:S^1_2\lhd U$, then there is a theory $\Th\supseteq U$ such that in every model of $\Th$, $\TB_N$
		 is definable with an atomic formula but $\Th$ does not syntactically define $\TB_N$ modulo $\Lang_U$.
	\end{stwierdzenie}
	\begin{proof}
		Fix $U$ and an interpretation $N: S^1_2\lhd U$.
		 Consider the extension of $U$ in the language with two additional predicates $T_1(x)$, $T_2(x)$ with all sentences of the form
		\[T_1(\num{\qcr{\phi}})\equiv \phi \vee T_2(\num{\qcr{\psi}})\equiv \psi,\]
		for $\phi,\psi\in\Lang$. Call this theory $2\TB_N$. Observe that if $(\M,T_1,T_2)\models 2\TB_N$, then 
		\[(\M, T_1)\models \TB_N \textnormal{ or } (\M,T_2)\models \TB_N.\]
		So it is sufficient to show that $2\TB_N$ does not define $\TB_N$. We shall show that $2\TB_N$ does not impose elementary equivalence. Fix any $n$ and let $\M_1\models U$ be an $\omega$-saturated model of $U$. In particular, $\M_1$ codes (in the sense of the defined membership predicate $\in_U$ satisfying AS) every consistent and complete extension of $U$. Assume that $k$ is the complexity of $U$-definition of $x\in_U y$. Let $\theta$ be a $(k+n,0)$-flexible formula over $U$. Without loss of generality assume that $\M_1\models \exists x \theta(x)$. By Proposition \ref{prop_flexible} there is $\M_2\succeq_{k+n} \M_1$ such that $\M_2\models \forall x \neg\theta(x)$. Let $c_1, c_2\in M$ be the codes of the theories of $\M_1$ and $\M_2$ respectively.
		Now put 
		\begin{align*}
			T_1^{\M_i}&:= \set{a\in M_i}{\M_i\models a\in_U c_1}\\
			T_2^{\M_i}&:= \set{a\in M_i}{\M_i\models a\in_U c_2}
		\end{align*}
		It follows that $(\M_i,T_1^{\M_i}, T_2^{\M_i})\models 2\TB_N$ and 
		\[(\M_1,T_1^{\M_1}, T_2^{\M_1})\preceq_n (\M_2, T_1^{\M_2},T_2^{\M_2}).\]
		However $\Th(\M_1)\neq \Th(\M_2)$, which ends the proof.
	\end{proof}
In fact, the above proof gives us a stronger statement:

 \begin{stwierdzenie}\label{prop_ogr_sem_tb}
		If $U$ is any sequential theory and $N: S^1_2\lhd U$, then there is a theory $\Th\supseteq U$ such that every complete extension of $\Th$ syntactically defines
		$\TB_N$ modulo $\Lang_U$ with an atomic formula but $\Th$ does not syntactically define $\TB_N$ modulo $\Lang_U$.
\end{stwierdzenie}

Corollary below answers the question of Albert Visser from \cite{viss19enayat}[Question 5.7]

	\begin{wniosek}\label{cor_nonres_tb}
		For every sequential theory $U$ and interpretation $N: S^1_2\lhd U$, $\TB_N$ does not have a restricted axiomatization modulo $U$.
	\end{wniosek}
	\begin{proof}
		Suppose there is an $n$-restricted theory $V$ such that $U+V\equiv \TB_N$. Let $\M_1$, $\M_2$,  $T_1^{\M_1}$, $T_1^{\M_2}$ be constructed as in the proof of Proposition \ref{prop_ogr_sem_tb_2}.
		 Since $(\M_1, T_1^{\M_1})\models \TB_N$ it follows that $(\M_1, T_1^{\M_1})\models V$ and since $V$ is $n$-restricted, then $(\M_2, T_1^{\M_2})\models V$.  Hence $(\M_2, T_1^{\M_2})\models \TB_N$, which is impossible, because $\num{\qcr{\exists x \theta(x)}}\in T_1^{\M_2}$ but $\M_2\models \forall x \neg\theta(x)$. 
	\end{proof}
	
	An observant reader might wonder whether there is anything interesting to say about the uniform version of a theory $\DEF_K$. One way to formulate this theory is via the following collection of sentences
    \[\forall x\forall y \bigl(D(\num{\qcr{\phi(v,w)}}, x,y) \equiv \exists !v\phi(v,y) \wedge \phi(x,y)\bigr).\]
    The intuitive meaning of $D(\num{\qcr{\phi(v,w)}}, a,b)$ is that $a$ is definable via the formula $\phi(x,b)$, where $b$ is a parameter. Let us call this theory $\textnormal{UDEF}_K$

    Actually, it is very easy to see that for every $U$, $\textnormal{UDEF}_K$ and $\USB_K$ are mutually definable. Indeed, the definition of $D(x,y,z)$ in $\USB_K$ is analogous to the case of $\DEF_K$ and $\USB_K$, and the definition of $S(\phi(v),y)$ in $\UDEF_K$
     is given by the formula:
    \[D(\qcr{v=w \wedge \phi(v)}, y,y).\]

	\subsection{Skolem functions}\label{subsect_skolem}
	
	We devote this section to the concept of Skolem functions. We treat it separately, since we do not have a clear semantical property which characterizes it up to the definability. However it is worth mentioning since it is classical and nicely fits into our hierarchy. As previously we work over a fixed sequential theory $U$ in a language $\Lang$ and $K$ is a distinguished interpretation $S^1_2\lhd U$.
	
	\paragraph{Skolem functions} S$_K(\Lang)$ extends $U$ with all sentences of the form
	\[\exists x \phi(x)\rightarrow \phi(H(\num{\qcr{\phi}})),\]
	where $\phi(x)$ is an $\Lang$-formula, with a single free variable and $H(x)$ is a fresh function symbol.
	
	\paragraph{Uniform Skolem functions} US$_K(\Lang)$ extends $\EA$ with all sentences of the form
	\[\exists x \phi(x,y)\rightarrow \phi(H(\num{\qcr{\phi}}, y), y),\]
	where $\phi(x,y)$ is an $\Lang$-formula with at most $x,y$-free and $H(x,y)$ is a fresh function symbol.
	
	\begin{stwierdzenie}
		For every sequential theory $U$ and every $K: S^1_2\lhd U$, $\TB_K(\Lang_U)$ and $\DEF_K(\Lang_U)$ are syntactically definable in $\textnormal{S}_K(\Lang_U)$ modulo $\Lang_U$. Similarly, $\USB_K(\Lang_U)$ is syntactically definable in $\textnormal{US}_K(\Lang_U)$ modulo $\Lang_U$.  
	\end{stwierdzenie}
	\begin{proof}
		We do the case of $\USB_K(\Lang_U)$, the rest of cases is very similar. Define
		\[S(\phi,x):= \num{1}\bigl(H(\theta_{\phi},x)\bigr),\]
		where $\theta_{\phi}(z,x):= \bigl((\num{1}(z)\wedge \phi(x))\vee (\num{0}(z) \wedge \neg\phi(x))\bigr)$ (we recall that $\num{n}(x)$ is the analogue of $x = \num{n}$ in the context of a general sequential theory).
		
	\end{proof}
	
	Now, we turn to separations between theories of truth, definability, satisfaction on one side and skolem functions on the other. We start with the comparison od truth and witnessing functions. The following lemma is a variant of Lemma \ref{lem_hierarchical_model}. We stress, that unlike its predecessor it holds without the additional assumption on the existence of infinitely many $\Sigma_n^*(\Lang)$ definable elements.
	
	\begin{lemat}\label{lem_hierarchical_model_2}
		For every r.e. sequential theory $U$ in  a finite language, every model $\M\models U$ and every $n\in \omega$, there is a model $\M\preceq_n \N$ such that $N = \bigcup_{i\in\omega}M_i$ and for every $l,m\in\omega$, there is $\theta(x)\in\Lang$ such that $\N\models \exists x\theta(x)$ but for every $\Sigma^*_l(\Lang_{M_m})$-definable element $a\in N$, $\N\models \neg\theta(a)$.
	\end{lemat}
	\begin{proof}
		The strategy is the same as in the proof of Lemma \ref{lem_hierarchical_model}. We build sequences $\was{\M_i}$, $\was{\theta_i}$ and $\was{n_i}$ which \rem{satisfy conditions $1.$ and $2.$ from the proof}
		and we switch condition $3.$ to the following one
		\begin{itemize}
			\item[3'.] for each $i$, $\M_{i+1}\models \exists x \theta_{i+1}(x)$ but for every $\Sigma_{n_i}^*(\Lang_{\M_i})$ definition $\phi(x)$, $\M_{i+1}\models \neg(\theta_{i+1}(x)\wedge \phi(x))$. 
		\end{itemize}
		It is a routine exercise to verify that for so defined sequence $\was{\M_i}_{i\in\omega}$, $\N:= \bigcup_{i\in\omega}\M_i$ satisfy the thesis of the lemma. Let us now show how to construct the relevant sequences. We define $\M_0 = \M, n_0 = n+1$
		 and choose $\theta_0$ to be any formula. Since $U$ is sequential,
		 there is $k\in\omega$ and $\was{\phi_i(x)}_{i\in\omega}$ such that for each $j$, $\phi_j$ is of complexity $\Sigma^*_k$, $\M\models \exists x \phi_j(x)$ and for any $i\neq j$, $\M\models \neg \bigl(\phi_i(x)\wedge \phi_j(x)\bigr)$ (for example $\phi_i$ defines the $E$-equivalence class corresponding to the natural number $i$). Assume that we are given $\M_i, n_i$. We let $\theta_{i+1}(x)$ be the $(n_i+1, k)$-flexible formula. Consider the following theory
		\[\Th:= U + (n_i+1)\mhyphen\Diag(\M_i) + \was{\theta_{i+1}(d)} + \set{\neg\theta_{i+1}(a)}{a\in M_i}.\]
		We claim that $\Th$ is consistent. For this, fix a finite fragment $A\subseteq \Th$ and  assume $a_1,\ldots, a_l$ are all elements from $\M_i$ such that $\neg\theta(a)$ occurs in $A$. Let $m$ be any natural number such that for all $i\leq l$, $\M_i\models \neg\phi_m(a_i)$. By the flexibility of $\theta_{i+1}$ there is a model $\M'\models U + \forall x (\theta_{i+1}(x)\equiv \phi_m(x))$ such that $\M_{i+1}\preceq_{n_i+1}\M'$. Clearly, $\M'\models A$. Now we set $\M_{i+1}$ to be any model of $\Th$ and put $n_{i+1}:= \Sigma^*(\theta_{i+1}) + 1$.
		 Observe that, by elementarity, every $\Sigma_{n_i}^*(\Lang_{M_i})$ definable element in $\M_{i+1}$ belongs to $M_i$. Hence, no witness for $\exists x \theta_{i+1}(x)$ is $\Sigma^*_{n_i}(\Lang_{M_i})$-definable in $\M_{i+1}$.
	\end{proof}
	
	\begin{tw}
		For every r.e. sequential $U$ in a finite language and every $K: S^1_2\lhd U$, $\TB_K(\Lang_U)$ does not semantically define $S_K(\Lang_U)$. 
	\end{tw}
	\begin{proof}
		Fix any model $\M\models U$ which is $\omega$-saturated. Let $\N$ be a model satisfying the thesis of Lemma \ref{lem_hierarchical_model_2}
		for $n$ equal to the complexity of the definition of $x\in y$ in $U$. In particular $\N$ is an $n-$elementary extension of $\M$. It follows that $\N$ codes $\Th(\N)$. Let $c\in M$ be such a code. Consequently $(\N,c)\models\TB_K[x \in_{\N} c/T(x)]$.
		 We claim that $\ElDiag(\N, c)$ does not syntactically define  $S_K$ modulo $\Lang_U$. Assume the contrary and let $k$ be the complexity of the definition of $H(x,y)$ (in a relational form).  Assume further that all the parameters used in the definition belong to $M_m$. Then $H(x,y)$ is $\Sigma_{k+n}^*(\Lang_{M_m})$-definable in $\N$.
		It follows that each non-empty parameter-free definable set $A\subseteq N$ contains a witness which is $\Sigma_{k+n}^*(\Lang_{M_m})$-definable in $\N$. This contradicts the property of $\N$.
	\end{proof}
	
	\begin{stwierdzenie}
		If $U\subseteq\Lang_{\PA}$ extends $\PA$, then $\DEF_{id}(\Lang_{\PA})$ defines $\textnormal{S}_{id}(\Lang_{\PA})$ and $\USB_{id}(\Lang_{\PA})$ syntactically defines $\textnormal{US}_{id}(\Lang_{\PA})$ modulo $\Lang_{\PA}$.
	\end{stwierdzenie}
	\begin{proof}
		Essentially this follows, since $\PA$ has definable Skolem functions. Let us do the case of $\textnormal{US}_{id}(\Lang_{\PA})$. Working in $\PA+\USB_{id}(\Lang_{\PA})$ we define
		\[H(\phi(x,v),y) = z := S(\qcr{\phi(x,v)\wedge \forall w<x\neg\phi(w,v)}, \pair{z}{y}).\] 
	\end{proof}
	
	As the next two theorems show, in the above result $\PA$ cannot be replaced by any of its finite fragments. Unfortunately we do not know exactly over which sequential theories, $\DEF_N(\Lang)$ semantically/syntactically  defines $\textnormal{S}_N(\Lang)$. Similarly, we do not know over which sequential theories, $\USB_N{\Lang}$ semantically/syntactically defines $\textnormal{US}_N(\Lang)$.
	However we can isolate a fairly natural class of theories over which the definability fails.
	
	\begin{tw}\label{stw_nondef_H_usb}
		For every $n\in\omega$, and every $U\subseteq \Lang_{\PA}$ consistent with $B\Sigma_n + \neg I\Sigma_n + \exp$, $\USB_{id}(\Lang_{\PA})$ does not parameter-free semantically define $\textnormal{S}_{id}(\Lang_{\PA})$ modulo $\Lang_{\PA}$.  
	\end{tw}
	\begin{proof}
		Fix $U$ and $n\in\omega$ and let $\M\models U$ be countable and recursively saturated. We claim that there is a (non-empty) parameter-free definable set $A\subseteq M$ such that $A$ contains no parameter-free definable element. Assume the contrary. Then, by the Tarski-Vaught test, $K(\M)\preceq \M$. In particular $K(\M)\models B\Sigma_n + \neg I\Sigma_n + \exp$ and every element of $K(\M)$ is parameter-free definable in $K(\M)$. Since $K(\M)$ has no proper elementary submodel, this literally contradicts \cite[Theorem 7.4.4]{kossakschmerl}. 
		
		So assume $A\subseteq M$ is a parameter-free definable subset of $\M$ with no parameter-free definable element. Since $\M$ is recursively saturated, by an easy back-and-forth argument,
		for each $a\in A$ there is an automorphism of $\M$ which moves $a$. Let $S:= \set{\pair{\qcr{\phi(x)}}{a}\in \omega\times M}{\M\models \phi(a)}$ be the Tarskian satisfaction class for $\M$. It follows that $(\M,S)\models \USB_{id}$. Moreover, the automorphism groups of $\M$ and $(\M,S)$ are the same (because each isomorphism preserve definable sets and pointwise fixes $\omega$). Suppose $\textnormal{S}_{id}(\Lang_{\PA})$ is parameter-free definable in $(\M,S)$. It follows that there is an element of $A$ which is parameter-free definable in $(\M,S)$. However, by the choice of $A$, $S$ and $\M$, every element of $A$ can be moved by an automorphism of $(\M,S)$. This concludes the proof.
	\end{proof}
	
	\begin{tw}\label{tw_nondef_US_in_UTB}
		For every $n\in\omega$, and every $U$ consistent with $B\Sigma_n + \neg I\Sigma_n + \exp$, $\USB_{id}(\Lang_{\PA})$ does not semantically define $\textnormal{US}_{id}$ modulo $U$.  
	\end{tw}
	\begin{proof}
		We modify the above proof of Proposition \ref{stw_nondef_H_usb}. We fix $U, n,\M$ exactly as above and claim that for each $a\in M$ there is an $\Lang_{\PA}\cup \was{a}$ definable subset of $M$ with no $\Lang_{\PA}\cup\was{a}$ definable element. The idea of the argument is as previously, only instead of $K(\M)$ we are using $K(\M,a)$ (the model consisting of $(\Lang_{\PA}\cup \was{a})$-definable elements of $\M$) and instead of \cite[Theorem 7.4.4]{kossakschmerl} we use \cite[Corollary 3.5]{kaye_model_char_PA}.
		
		Now we let $S$ be defined exactly in the same way as previously. We fix $a\in M$ and show that there is no definition of $H$ which uses $a$ as a parameter. Assume that there is such a definition. Let $A$ be the $a$-definable set with no $a$-definable element. It follows that for every element $b\in A$ there is an automorphism $f\in \Aut(\M)$ which fixes $a$ and moves $b$. As previously it follows that no element in $A$ is $a$-definable in $(\M,S)$, which ends the proof. 
	\end{proof}
	
	\begin{stwierdzenie} \label{stw_S_nie_definiuje_USB}
		For every sequential theory $U$ and $K: S^1_2\lhd U$, $\textnormal{S}_K(\Lang_U)$ does not semantically define $\USB_K(\Lang_U)$ modulo $\Lang_U$.
	\end{stwierdzenie}
	\begin{proof}
		This follows from the fact that each sequential theory has a model $\M$ such that for some $c\in M$, $(\M,c)\models \textnormal{S}_N(\Lang_U)[\tuple{x,y} \in_{\M} c / H(x) = y ]$.
		 Indeed, let $\M\models U$ be any recursively saturated model of $U$. Consider the following recursive type $p(x)$ (recall (Convention \ref{konwencja_arytmetyka_w sekwencji}) that we treat $\num{\qcr{\phi(y)}}$ as a predicate
		which isolates $E$-equivalence class of elements corresponding to the numeral naming $\qcr{\phi}$, so $\num{\qcr{\phi(y)}}(z)$ means that $z$ is equal to the numeral naming $\qcr{\phi(y)}$):
		\[\set{\exists y \phi(y)\rightarrow \exists !z\exists! y (\phi(y)\wedge \num{\qcr{\phi(y)}}(z)\wedge \pair{z}{y}\in_{\M} x)}{\phi(y)\in\Lang}.\]
		By the sequentiality of $U$, $p(x)$ is a type. Any elements realizing this type codes a witnessing function for $\Lang$-formulae. However, a predicate satisfying $\USB$ can never be definable with a parameter. 
	\end{proof}

	\subsection{$\PA$ case}\label{subsect_PA}

	In this subsection and until the end of the paper (i.e. also in Section \ref{sect_truth_rec_sat}), we specialize to languages which extend $\Lang_{\PA}$ with at most countably many new predicate symbols and theories that extend $\PA$ and prove full scheme of induction for their respective languages. In particular, we will omit the mention of the interpretation of $S^1_2$ in our theory, as it will be assumed to be identity. 
	As it will be important to distinguish between the situation in which we have full induction for the truth predicate and the one in which induction is assumed only for arithmetical formulae, we shall employ the following notation: the minus sign in the superscript will indicate that we do not extend induction axioms to the whole language.
	 In particular, $\USB^-(\Lang)$ 
	 denotes simply $\PA + \USB_{id}(\Lang)$ (in the terminology of Section \ref{sect_prelim_truth}) and $\USB(\Lang)$ is its extension with induction axioms for the extended language. The collection of all induction axioms for a language $\Lang$ will be denoted $\Ind_{\Lang}$.

	We start with a proposition that is properly weaker than already proved Theorem \ref{tw_nonuniform_impos_elem_utb}, but we give it here because in this context it has a particularly illustrative proof.
	
	\begin{stwierdzenie}
		If $\PA\subseteq U$ is an r.e. theory in a finite language which uniformly imposes $\Lang$-elementarity (elementary equivalence, preserves definability) and such that $U\vdash \Ind_{\Lang_U}$, then $\USB(\Lang)$ ($\TB(\Lang)$, $\DEF(\Lang)$) is syntactically definable in $U$.
	\end{stwierdzenie}
	\begin{proof}
		We sketch the proof for $\USB(\Lang)$ (the proof easily adapts to the other two cases). We use the following convention: if $\phi(x)$ is an arbitrary formula, then $\Con(\phi(x))$ denotes the assertion that the theory consisting of all those axioms that satisfy $\phi(x)$ is consistent. Choose $n$ witnessing that $U$ uniformly imposes $\Lang$-elementarity. 
		
		
		Consider the following formalized theory ($\True_k(x)$ is a partial truth predicate for $\Sigma^*_k$ sentences, which is based on $\Sat_k(x.y)$):
		\[\alpha(x):= \True_{n+1}(x) \vee \bigl(U(x)\wedge \Con(U\restr{x}+\True_{n+2})\bigr).\]
		In the above $U(x)$ is an elementary predicate which represents the axiom set of $U$ and $U\restr{x}$ denotes the theory consisting of all axioms of $U$ of G\"odel code $\leq x$ and $\Con(U\restr{x}+\True_{n+2})$ expresses that every $\Sigma^*_{n+2}$ sentence provable in $U\restr{x}$ is true.  Since $U$ contains full $\PA$, $U$ proves cut-elimination theorem for first-order logic.
		In particular, for every $k,n$ separately, by considering cut-free proofs, inside $U$ we can reprove the standard inductive argument and show that $U\restr{k} + \True_{n+1}$ is consistent (i.e. axioms of $U$ smaller than $k$ are consistent with all $\Sigma^*_{n+1}$ true sentences). Hence $U\vdash\Con(\alpha(x))$ (i.e. the set of axioms defined by $\alpha(x)$ is a consistent theory.)
		and, for every $k$ separately, 
		\begin{equation}\label{equat_N_models_U}
			U\vdash U\restr{k}\subseteq \alpha.
		\end{equation}
		By the Arithmetized Completeness Theorem (see \cite[Section 13.2]{kaye}) there is a definable model $\N = (N, +^{\N}, \cdot^{\N}, 0^{\N}, 1^{\N})$ which comes with a definable truth predicate $\models^{\N}$ such that $U$ proves
		\begin{equation}\label{equat_N_models_alpha}
			\forall \phi\bigl(\alpha(\phi)\rightarrow \N\models^{\N} \phi\bigr)
		\end{equation}
		Let $\iota$ denote the canonical initial embedding of the universe onto an initial segment of $\N$, defined recursively $\iota(0) = 0^{\N}$, $\iota(x+1) = \iota(x)+^{\N} 1^{\N}$ 
		(i.e. $x$ is mapped to the value if the $x$-th numeral as computed in $\N$) and define
		\[T_U(\phi, x):= \N\models^{\N}\phi[\iota(x)].\]
		We check that $U\vdash \USB^-(\Lang)[T_U/T]$. Pick any model $\M\models U$. Then the definable model $\N$ gives rise to a true model (which will be denoted with the same letter) and $\models^{\N}$ coincides with the usual satisfaction class for $\N$. Hence, by \eqref{equat_N_models_U} and \eqref{equat_N_models_alpha}, $\N\models U$. By the properties of $\iota$, $\iota[\M]\subseteq_e \N$. Since $\M\models \True_{n+1}\subseteq \alpha$, $\iota[\M]\preceq_n \N$. It follows that $\iota[\M\restr{\Lang}]\preceq \N\restr{\Lang}$, because $U$ imposes $\Lang$-elementarity. In particular, for an arbitrary $a\in M$ and $\phi(x)\in\Lang$ the following are equivalent
		\begin{enumerate}
			\item $\M\models T_U(\phi,a)$
			\item $\M\models \N\models^{\N}\phi[\iota(a)]$
			\item $\N\models \phi[\iota(a)]$
			\item $\iota[\M]\models \phi[\iota(a)]$
			\item $\M\models \phi(a)$
		\end{enumerate}
		The equivalence between $2.$ and $3.$ is the absoluteness of the satisfaction relation and the equivalence between $3.$ and $4.$ follows, since $\iota[\M\restr{\Lang}]\preceq \N\restr{\Lang}$.
	
	\end{proof}
	
	In the separations below we shall use the following observation:
	
	\begin{stwierdzenie}\label{stw_def_bounded}
		Let $\Lang$ be at most countable language. Suppose $(\M,D)\models \PA + \DEF(\Lang)$ is nonstandard. Then the $\Lang$-definable elements are bounded in $\M$.	
	\end{stwierdzenie}
	\begin{proof}
		We use overspill on the formula
		\[\psi(x):= \exists c \forall \phi(z)<x \forall y \bigl(D(\phi(z), y)\rightarrow y<c\bigr).\]
		Clearly $(\M,D)\models \psi(n)$ for every $n\in\omega$. Hence $\psi(d)$ holds for some $d\in M\setminus\omega$ and the witness for $\psi(d)$ bounds all the definable elements.
	\end{proof}
	
	\begin{stwierdzenie}
		If $U\supseteq \PA$ is an r.e. theory in at most countable language $\Lang_U$ such that $U\vdash \Ind_{\Lang_U}$, then $U+\USB(\Lang_U)\nleq^s_{\Lang_U} U+\DEF(\Lang_U)\nleq^s_{\Lang_U} U+\TB(\Lang_U)$.
	\end{stwierdzenie}
	\begin{proof}
		The second separation follows, since there is a model $(\M,T)\models U+\TB(\Lang_U)$ in which the $\Lang_U$-definable elements are cofinal. This model cannot be expanded to a model of $\DEF(\Lang_U)$, by Proposition \ref{stw_def_bounded}. To get this model fix any $(\M, T)\models \TB(\Lang_U)$, which is not a model of true arithmetic. Let $c\in SSy(\M)$ be the code of the theory of $\M$ (this element exists because we have induction for the truth predicate). Then $\mathcal{K}(\M)$ (the submodel consisting of $\Lang_U$-definable elements) is nonstandard and elementary in $\M$. By the Gaifman splitting theorem for $U$ (se \cite[Corollary 7.10]{kaye}) there is $\N$ such that 
		\[\mathcal{K}(\M)\preceq_{cof}\N\preceq_e \M.\]
		(in fact, $\N$ is simply the supremum of $\mathcal{K}(\M)$ in $\M$). Since $\N\preceq_e\M$, we may assume that $c\in N$, $(\N, c)\models \TB(\Lang_U)[x \in c/T(x) ]$.
		 However in $\N$ the definable elements are cofinal, because $\mathcal{K}(\M)\preceq_{cof}\N$.
		
		The first separation follows, since there is a model of $\DEF(\Lang_U)$ whose $\Lang_U$-reduct is not short recursively saturated, and $\USB(\Lang_U)$ imposes $\Lang_U$-recursive saturation (see e.g. \cite{wcislyk_models_weak}). To see that such a model exists, fix any nonstandard $(\M,D)\models \DEF(\Lang_U)$ and let $d$ be such that for every $\phi(x)$ and every $a\in M$
		\[(\M,D)\models \pair{\num{\qcr{\phi(x)}}}{a}\in d \equiv D(\num{\qcr{\phi(x)}}, a).\]
		The existence of $d$ follows once again by overspill. Take $\N:= \mathcal{K}(\M,d)$ (i.e. the submodel of $\M$ consisting of elements $(\Lang_U \cup\was{d})$-definable with parameter $d$). Then $d\in\N\preceq_{\Lang_U}\M$ and $\N$ is not short recursively saturated, since it does not realise the type $p(x,d)$

		\[\set{x<d \wedge \bigl(\exists!y\phi(y,d)\rightarrow \neg \phi(x,d)\bigr)}{\phi(y)\in\Lang_{U}}.\]
		However, $\N$ expands to a model of $\DEF(\Lang_U)$, because in $\N$ the predicate satisfying $\DEF(\Lang_U)$ is definable with a parameter via the formula
		\[D(x,y):= \Form^{1}(x) \wedge \pair{x}{y}\in d.\]
		where $\Form^1(x)$ says that $x$ is a formula with one free variable.
	\end{proof}	
	
	Recall that $\leq^s_{\Lang}$ and $\leq^m_{\Lang}$ denote syntactical and semantical definability modulo $\Lang$ (respectively).
	\begin{wniosek}
		$\TB\leq^s_{\Lang_{\PA}}\DEF\leq^s_{\Lang_{\PA}}\USB$ but $\TB\ngeq^m_{\Lang_{\PA}} \DEF\ngeq^{m}_{\Lang_{\PA}} \USB$.
	\end{wniosek}

\begin{uwaga} \label{rem_separation_noninductive_PA}
	As a corollary we also obtain the above separations and reductions also for noninductive version of the above theories. For the separations: notice that if there were reductions in the noninductive case, they would simply imply the respective reductions in the fully inductive case. On the other, the reductions which we obtained for the inductive versions of the theories, actually only used the axioms of $\DEF^-, \TB^-$ and $\USB^-$, respectively, not the induction scheme for the predicates $D$, $T$, and $S$. 
\end{uwaga}
	
	\subsubsection{Definability in no model}\label{subsect_PA_def_inno}
	
	In this subsection we strengthen the non-definability results obtained above. Once again, we switch to non-inductive variants of $\TB$, $\DEF$ and $\USB$, which will be denoted simply $\TB^-$, $\DEF^-$ and $\USB^-$, respectively.

    In the proposition below we are going to use the following result due to Stuart Smith (\cite[Theorem 3.11]{Smith_extendible}, and the remark immediately after the proof on p. 351)

    \begin{tw}\label{Smith_omega}
    Suppose that $\M\models \PA$ is nonstandard. Then for every $\N$ such that $\M\preceq \N$ and every $A\subseteq \omega$, $(\M,A)\preceq (\N,A)$.
    \end{tw}
    In the above, $(\M,A)$ denotes the expansion of a structure $\M$ to a language with a fresh predicate interpreted as set $A$.
	\begin{stwierdzenie}\label{tw_ost_niedef_TB_UTB}
		In every nonstandard model
		$\M\models \PA$ there is a set $T$ such that $(\M,T)\models\TB^-(\Lang_{\PA})$ but $\ElDiag(\M,T)$ does not syntactically define $\USB^-(\Lang_{\PA})$.
	\end{stwierdzenie}
	\begin{proof}
		Fix $\M\models \PA$ and let $T = \set{\qcr{\phi}\in\omega}{\M\models \phi}$. We claim that $\ElDiag(\M,T)$ does not define syntactically define $\USB^-(\Lang_{\PA})$ modulo $\Lang_{\PA}$. Assume the contrary and pick any formula $S(x,y)\in \Lang_M$. Let $a$ be a parameter used in this definition and without loss of generality assume that $a$ is nonstandard. 
		Let $\N$ be any proper elementary end extension of $\M$ (which exists by the Macdowell-Specker Theorem \cite[Theorem 8.6]{kaye}) and let $\mathcal{K}(a)$ be the submodel of $\M$ generated by the elements which are definable in $\M$ from $a$. Moreover, let $\M'$ be the submodel of $\M$ with the domain $\set{d\in M}{\exists x\in \mathcal{K}(a) \M\models d<x}$.
		By the result of Gaifman (see \cite[Corollary 7.10]{kaye}) $\M'\preceq_e \M$ (we do not assume that this extension is proper.)
		 We have 
		\[\M'\preceq_e \M\precneqq_e \N,\]
		so by Theorem \ref{Smith_omega}
		for $A=T\subset\omega$
		\[(\M',T)\preceq_e(\M,T)\precneqq_e (\N,T).\]
		Note that $\omega$ is definable in all the three models with a formula $\phi(x):= \exists y \bigl(x\leq y \wedge T(y)\bigr)$. Moreover $S(x,y)$ defines $\USB^-(\Lang_{\PA})$ in both $(\M',T)$ and $(\N,T)$ and hence 
		\begin{equation}\label{elementarity}\tag{E}
			(\M',S^{\M'}, \omega)\precneqq_e(\N,S^{\N},\omega)
		\end{equation}
		(recall that for a formula $\phi$ and a model $\mathcal{U}$, $\phi^{\mathcal{U}}$ denotes the set of elements which satisfy $\phi$ in $\mathcal{U}$.) However, the $\Lang_{M'}\cup \was{T}$ formula
		\[\zeta(x,a):= \exists \psi(v_0,v_1)\in \omega\bigl[ S(\qcr{\psi(v_0,\num{a}/v_1)}, x)\wedge \forall z\bigl(S(\qcr{\psi(v_0, \num{a}/v_1)}, z)\rightarrow x=z\bigr)\bigr]\]
		in both $(\M',S^{\M'},\omega)$ and $(\N, S^{\N}, \omega)$ defines the set of $\Lang_{\was{a}}$-definable elements of $\M'$ and $\N$, respectively. Since $(\M',S^{\M'},\omega)\models \forall y\exists x>y \zeta(x,a)$ and $(\N,S^{\N}, \omega)\models \exists y \forall x>y \neg \zeta(x,a)$, this contradicts \eqref{elementarity}.
	\end{proof}
	
	\begin{stwierdzenie}
		For every model $\M\models \PA$ the following are equivalent:
		\begin{enumerate}
			\item $\M\models \Th(\mathbb{N})$
			\item For every $T\subseteq M$, if $(\M,T)\models \TB^-(\Lang_{\PA})$, then $\ElDiag(\M,T)$ syntactically defines $\DEF^-(\Lang_{\PA})$.
		\end{enumerate}
	\end{stwierdzenie}
	\begin{proof}
		Fix any $\M\models \PA$. Assume first that $\M\models \Th(\mathbb{N})$ and take any $T\subseteq M$ such that $(\M,T)\models \TB^-(\Lang_{\PA})$. The following formula is then a definition of $\DEF^-(\Lang_{\PA})$
		\[D(\phi(x),y):= T(\qcr{\exists! x \phi(x)})\wedge T(\qcr{\phi(\dot{y})})\wedge \forall z< y \neg T(\qcr{\phi(\dot{z})}).\]
		We prove that it works. Assume first $\M\models \exists !x \phi(x)\wedge \phi(a)$. Since $\M\models \Th(\mathbb{N})$, $a\in \omega$. In particular $\phi(\num{a})$ is a standard sentence and by $\TB^-(\Lang_{\PA})$ $(\M,T)\models T(\exists!x\phi(x))\wedge T(\qcr{\phi(\num{a})})$. Since every $b<a$ is also a standard number, then for any such $b$ we have $(\M,T)\models \neg T(\qcr{\phi(\num{b})})$. Now assume $(\M,T)\models D(\phi,a)$. By $\TB^-(\Lang_{\PA})$ it follows that $\M\models \exists!x\phi(x)$. We claim that $\M\models \phi(a)$. Assume the contrary. Since $(\M,T)\models T(\qcr{\phi(\num{a})})\wedge \neg\phi(a)$ then $a$ is nonstandard. However, since $\M\models \Th(\mathbb{N})$ there is the unique $n\in\omega$ such that $\M\models \phi(\num{n})$. 
		It follows that $(\M,T)\models \exists x<a T(\qcr{\phi(\dot{x})})$, contradicting $D(\phi,a)$.
		
		Now assume that $\M$ fails to satisfy $\Th(\mathbb{N})$ and let $T:= \set{\phi\in\omega}{\M\models \phi}$. Let $\mathcal{K}$ be the submodel of $\M$ consisting of $\Lang_{\PA}$ (i.e. parameter-free) definable elements of $\M$. Let $\M'$ satisfy $\K\preceq_{cf}\M'\preceq_e \M$. Let $\N$ be any proper elementary end-extension of $\M$. Now mimic the proof of Theorem \ref{tw_ost_niedef_TB_UTB} to conclude that if $\ElDiag(\M,T)$ defines $\DEF^-(\Lang_{\PA})$,
		 then the extensions
		\[(\M', T)\subseteq_e (\M,T)\subsetneq_e (\N,T)\]
		cannot all be elementary, contradicting Theorem \ref{Smith_omega}.
	\end{proof}

	Now we comment on the relation between $\DEF^-(\Lang_{\PA})$ ans $\USB^-(\Lang_{\PA})$. Here, unlike in the previous two cases, we lack the full characterisation of the class of models whose every expansion to a model of $\DEF^-(\Lang_{\PA})$ defines $\USB^-(\Lang_{\PA})$.
	
	\begin{stwierdzenie}
		Assume $\M\models \PA$ is such that every element in $\M$ is definable without parameters. Let $D\subseteq M^2$ be such that $(\M,D)\models \DEF^-(\Lang_{\PA})$. Then $\USB^-(\Lang_{\PA})$ is syntactically definable in $\ElDiag(\Lang_{\PA})$.
	\end{stwierdzenie}
	\begin{proof}
		Let $\M\models \PA$ be prime and assume $(\M,D)\models \DEF^-(\Lang_{\PA})$. Let $T(x)$ be a definition of a truth predicate satisfying $\TB^-(\Lang_{\PA})$ in $(\M, D)$ (as in Proposition \ref{stw_latwa_definiowalnosc}). We define $S(\phi,x)$ in the following way
		\[\exists \psi \bigl(D(\psi,x)\wedge \forall y<\psi \neg D(y,x)\wedge T(\qcr{\forall x (\psi(x)\rightarrow \phi)})\bigr).\]
		The correctness of the definition is guaranteed by the fact that every element of $\M$ has a standard definition without parameters.
	\end{proof}

	Finally we notice that Theorem \ref{tw_ost_niedef_TB_UTB} fails for fragments of $\PA$ of restricted quantifier complexity. Let $\PA\restr{n}$ denote the set of consequences of $\PA$ of complexity $\Sigma_n$.
	
	\begin{stwierdzenie}
		For every $n$ there is a model $\M\models \PA\restr n$ such that for every $T\subseteq M$ such that $(\M,T)\models \TB^-(\Lang_{\PA})$, $\USB^-(\Lang_{\PA})$ is syntactically definable in $\ElDiag(\M,T)$.
	\end{stwierdzenie}
	\begin{proof}
		Fix $n$ and any model $\N\models \PA$. Let $\mathcal{K}^n$ be the submodel of $\M$ consisting of $\Sigma_n$-definable elements. Then, by \cite[Theorem 10.1]{kaye} $\K^n\preceq_n \M$, hence $\K^n\models \PA\restr{n}$. Moreover there is a formula $\phi(x)$ such that $\phi^{\K^n} = \omega$. Fix any $T$, a subset of the universe of $\K^n$
		 such that $(\K^n,T)\models \TB^-(\Lang_{\PA})$. Then it is easy to check that the following formula defines $\USB^-(\Lang_{\PA})$ in $(\K^n,T)$
		\[\theta(\psi(v),a):= \exists \chi(v) \in\Sigma_n \bigl(\phi(\chi) \wedge \Sat_{n}(\chi, a)\wedge \forall x\neq a\neg\Sat_n(\chi,x) \wedge T\bigl(\qcr{\forall v \bigl(\chi(v)\rightarrow \psi(v)\bigr)}\bigr).\]
	\end{proof}

	\section{Truth and recursive saturation}\label{sect_truth_rec_sat}
	
	We recall that in this section all theories in considerations are extensions of $\PA$ (possibly in an extended language) and $^-$ superscript signalizes that induction is assumed only for arithmetical formulae. The following definition will be handy in the particular context of $\PA$ as the base theory.
	
	\begin{definicja}
		In what follows
		the variables $s,t,u$ will range over G\"odel codes of arithmetical terms and variables $v$, $v_i$, $w$ - over codes of variables. For example the quantification $\forall t\psi$ should be understood as $\forall x \bigl(\ClTerm(x) \rightarrow \psi[x/t]\bigr)$. Similarly $\phi,\psi$ will range over codes of arithmetical sentences, $w,v,$ over codes of variables and $\phi(v)$, $\psi(v)$ over codes of arithmetical formulae with at most one free variable shown. For an arbitrary closed term $t$, $\val{t}$ denotes the value of $t$. The function $t\mapsto \val{t}$ is provably total in $\PA$, so we shall use the \rem{symbol $\val{x}$}
		as if it were a function symbol in our language. $\phi\in\dpt(x)$ is an $\Lang_{\PA}$ formula in variables $\phi$ and $x$ which expresses that $\phi$ is a $\Lang_{\PA}$ formula of syntactical depth at most $x$, i.e. each path in the syntactic tree of $\phi$ has length at most $x$ (compare \ref{def_complexity}).
		
		$\CT^-(x)$ denotes the conjunction of the axiom of $I\Delta_0+\exp$ and the following sentences of $\Lang_{\PA}\cup \was{T}$:
		\begin{enumerate}
			\item[CT1] $\forall s,t \bigl(T(\qcr{s=t})\equiv \val{s} = \val{t}\bigr)$ 
			\item[CT2] $\forall \phi\in\dpt(x)\forall\psi\in\dpt(x) \bigl(T(\qcr{\phi\vee\psi})\equiv T(\phi)\vee T(\psi)\bigr)$.
			\item[CT3] $\forall \phi\in\dpt(x)\bigl(T(\qcr{\neg\phi})\equiv \neg T(\phi)\bigr)$. 
			\item[CT4] $\forall v\forall \phi(v)\in\dpt(x)\bigl(T(\exists v \phi)\equiv \exists y T(\qcr{\phi(\dot{y}/v)})$.
		\end{enumerate}
		In the above, as introduced in Convention \ref{konwencja_arytmetyka_w sekwencji}, $\qcr{\neg\phi}$ denotes a definable function in variable $\phi$ that given a sentence $\phi$ returns the G\"odel code of $\neg\phi$. The use of $\qcr{\cdot}$ in axioms $\CT1$, $\CT2$ and $\CT4$ should be understood analogously.
		
		$\CT^-$ extends $\PA$ with the sentence $\forall x \CT^-(x)$ and $\CT(x)$ is the theory extending $\PA$ with $\CT^-(x)$ and full induction scheme for $\Lang_{\PA}\cup \was{T}$.
	\end{definicja}

	\begin{uwaga}
		For theories $U$ extending $\PA$ (in fact much less is needed), $\UTB^-_{id}(\Lang_{U})$ is mutually syntactically definable modulo $\Lang_{\PA}$ with the following term variant of $\UTB^-_{id}(\Lang_{U})$,
		 consisting of the following biconditionals:
		\[\forall t \bigl(T(\qcr{\phi(t)})\equiv \phi(\val{t})\bigr).\]
		In the above, $\phi(x)$ is an arbitrary formula of $\Lang_{\PA}$ and $\qcr{\phi(t)}$ denotes the effect of formal substitution of $t$ for $x$ in $\phi(x)$. The definability holds, because sufficiently strong arithmetical theories prove that every object from the universe can be named by a closed term. The above will be our official definition of $\UTB_{id}(\Lang_{U})$ in this section. 	 Arguably, it is the best known version of the considered theories of satisfaction. Following the conventions from the previous section, we will omit the subscript indicating the identity interpretation of $S^1_2$ in our theories. Moreover, we will sometimes also omit the mention of the arithmetical language $\Lang_{\PA}$.
	
	\end{uwaga}
	
	It is a well known fact that if a model $(M,T)$ satisfies $\UTB$, then its arithmetical part $M$ is recursively saturated. Surprisingly, as shown by \cite{lachlan}, the same holds if we assume that $(M,T)$ satisfies $\CT^-$, a theory with compositional axioms for arithmetical sentences and no induction.
	
	Roman Kossak, in \cite{kossak_cztery_problematy}, has shown that a partial reverse holds. The proof makes a crucial use of the MacDowell-Specker Theorem (see \cite{kossakschmerl}). Recall that a model $\N$ is a conservative extension of $\M$ if $\M$ is a submodel of $\N$ and for every $A\subseteq N$, if $A$ is definable in $\N$ with parameters, then $A\cap M$ is definable in $\M$ with parameters.

    \begin{tw}[Macdowell-Specker]
    Let $\Th$ be a theory in at most countable language $\Lang$ such that $\Th\vdash \PA$  and $\Th$ proves all the instantiations of the induction scheme with formulae of $\Lang$. Then every model of $\Th$ has a proper, conservative, elementary extension.
    \end{tw}
	
	\begin{tw}[Kossak] \label{tw_kossak_indukcja_vs_nasycenie}
		Let $U$ be a theory in a countable language, containing $\PA$ and the full induction scheme (for the extended language). Then, if $U$ imposes (short) $\Lang_{\PA}$-recursive saturation, then $U$ semantically defines $\UTB(\Lang_{\PA})$ modulo $\Lang_{\PA}$.
	\end{tw}
	
	\begin{proof}
		Let $\M \models U$ and let $\M \precneqq_e \M'$ be a proper, conservative and elementary end-extension of $\M$. $\M'$ exists since $U$ is a theory in countable language satisfying full induction, so the Macdowell-Specker Theorem holds for $U$. 
		
		Fix an arbitrary $a \in M' \setminus M$. Let $s$ be an arbitrary element which realizes the following recursive type:
		\[p(x):= \set{\forall t<a\bigl( \qcr{\phi(t)}\in x \equiv \phi(\val{t})\bigr)}{\phi\in \Lang_{\PA}}.\]
		In particular, $s$ codes the set of true arithmetical sentences of standard complexity with terms smaller than $a$. Let $T = \set{ x \in M }{\M' \models x \in s}.$ Then $T$ is a class in $\M$, hence it is definable in $\M$ by conservativeness of the extension. By definition, $(\M,T)$ satisfies $\UTB^-(\Lang_{\PA})$. Since $\M$ satisfies full induction, actually $(\M,T) \models \UTB(\Lang_{\PA})$.
	\end{proof}
	
	In the above argument, induction is used as a crucial ingredient in order to obtain a conservative end-extension. If a model of a fragment of $\PA$ does not satisfy full induction, it can never have an elementary end-extension at all, so the argument as-is breaks down completely. This motivates the following question: Does Theorem \ref{tw_kossak_indukcja_vs_nasycenie} hold if we drop the assumption that $U$ proves the induction scheme for its own language? Let us start with a simple example showing that the strict analogue cannot be true for rather general reasons. 
	
	\begin{uwaga} \label{rem_theory_of_realisers}
		Consider the following theory $\textnormal{RSAT}(\Lang_{\PA})$ in the language extending $\Lang_{\PA}$ with a fresh ternary predicate $R(p,x,y)$. $\RSAT(\Lang_{\PA})$ extends $\PA$  with the following axioms of optimality (OP) and nonemptiness (NE):
		\begin{enumerate}
			\item[OP] $\exists x  \bigwedge_{i = 1}^n \phi^p_i(y_1,\ldots, y_k,x) \rightarrow \forall x\left( R(\num{p}, x, \tuple{y_1,\ldots,y_k}) \rightarrow \bigwedge_{i=1}^n  \phi^p_i(y_1,\ldots,y_k,x)\right),$
			\item[NE] $\exists x R(\num{p}, x,\tuple{y_1,\ldots,y_k})$,
		\end{enumerate}
		where $\num{p}$ is the G\"odel code of a total p.r. function $i\mapsto \phi_i^p$ and $\tuple{y_1,\ldots,y_k}$ is an arithmetical code for a sequence $y_1,\ldots, y_k$. To improve readability we shall write $R_p(x,y)$ instead of $R(p,x,y)$.
		
		We think of $R_p$ as ``type-realisers.'' Each $x$ such that $R_p(x,\bar{y})$
		holds satisfies as large a portion of  the type $p$ as possible (with parametres $\bar{y}$).
		
		Clearly, any model of $\RSAT(\Lang_{\PA})$ is recursively saturated and, conversely, any recursively saturated model of $\PA$ expands to a model of $\RSAT(\Lang_{\PA})$. In particular, if $M$ is a recursively saturated rather classless model of $\PA$ (see \cite{kossakschmerl}), then it expands to a model of RSAT, but not to a model of $\UTB$. (Notice that in the definition of $\RSAT(\Lang_{\PA})$ we do not assume induction for the extended language).
	\end{uwaga}

	The obstruction discussed in the above remark is rather general in nature and shows that a direct strengthening of Theorem \ref{tw_kossak_indukcja_vs_nasycenie} is not possible. However, we might still hope that Kossak's result essentially holds in greater generality, in that the only way to impose recursive saturation in a theory is ``essentially'' to use a truth predicate. This intuition turns out to be correct, at least for theories extending a sufficiently strong arithmetic. The following strengthening of Kossak's theorem is true:
	
	\begin{tw} \label{th_recursive_saturation_implies_truth}
		Every theory in a countable language which extends $\PA$ and imposes $\Lang_{\PA}$-recursive saturation semantically defines $\UTB^-(\Lang_{\PA})$ modulo $\Lang_{\PA}$.
	\end{tw}
	Notice that if $U$ contains full induction and semantically defines $\UTB^-$($\Lang_{\PA})$ modulo $\Lang_{\PA}$, then in each model the definable truth predicate automatically satisfies full induction scheme, hence the above theorem is indeed a strengthening of the result by Kossak. 
	
	Before we present the proof of the result, let us note a counterexample to a stronger thesis, since it actually precedes Theorem \ref{th_recursive_saturation_implies_truth} and partly serves as an inspiration for its proof. We could ask whether we could strengthen the conclusion of the theorem to \emph{syntactical definability}. It turns out that such a strengthening does not hold.  
	
	\begin{tw} \label{th_theory_with_rec_sat_without_utb}
		There exists a theory $U \supseteq \PA +\Ind_{\Lang_U}$ which imposes $\Lang_{\PA}$-recursive saturation but does not syntactically define $\UTB^-(\Lang_{\PA})$ modulo $\Lang_{\PA}$.
	\end{tw}
	
	\begin{proof}
		Let $(\phi_j)_{j\in\omega}$ be a primitive recursive enumeration of all arithmetical formulae. Consider the following theory $U$ in the arithmetical language with additional predicates $T_i, T_{\omega}, i \in \omega$. $U$ contains $\PA$, full induction for the extended language, and the following axioms:
		\begin{displaymath}
			\neg \forall t \Big(T_j (\qcr{\phi_i(t)}) \equiv \phi_i(\val{t})\Big) 
			\rightarrow \forall t \Big(T_{\omega} (\qcr{\phi_j(t)}) \equiv \phi_j(\val{t})\Big),
		\end{displaymath}
		where $i, j \in \omega$.

		We claim that for every model $\M \models U$, the arithmetical reduct of $\M$ is recursively saturated. Indeed, in every such model, there exists a predicate $S$ satisfying full $\UTB(\Lang_{\PA})$, since if a predicate $T_j$ does not satisfy all the axioms of $\UTB(\Lang_{\PA})$, then $T_{\omega}$ satisfies the instance of the uniform biconditional scheme for the formula $\phi_j$. Hence either one of $T_j$, $j<\omega$, satisfies the axioms of $\UTB(\Lang_{\PA})$, or none of them does, in which case $T_{\omega}$ satisfies it.  
		
		We further claim that $U$ does not syntactically define $\UTB(\Lang_{\PA})$ modulo $\Lang_{\PA}$. Indeed, suppose that there exists such a definition $\phi$. Then $\phi$ uses only finitely many predicates $T_i$. Fix any $N \in \omega$ large enough so that all such predicates have index $i< N$ or $i=\omega$.
		
		Consider a model of $U$ such that the arithmetical part is standard (i.e., isomorphic to $(\mathbb{N},+, \times,0)$) and the predicates $T_i$ are interpreted as follows: all predicates $T_i$ for $i \leq N$ are interpreted as empty sets, $T_{\omega}$ is defined as the partial truth predicate for formulae $\phi_i, i \leq N$:
		\begin{equation*}
			T_{\omega} := \set{\qcr{\phi_i(t)} \in \Sent_{\Lang_{\PA}}(\mathbb{N})}{i \leq N \wedge \mathbb{N} \models \phi_i(\val{t})},
		\end{equation*}
		and all $T_j$ for $j > N$ are defined as the standard arithmetical truth predicate, i.e.:
		\begin{equation*}
			T_j := \set{\qcr{\phi(t)} \in \Sent_{\Lang_{\PA}}(\mathbb{N})}{\mathbb{N} \models \phi(\val{t})}.
		\end{equation*}
		One can check that the model $(\mathbb{N},T_i,T_{\omega})_{i < \omega}$ satisfies the axioms of the theory $U$. On the other hand, in that model the predicates $T_i, i \leq N$ and $T_{\omega}$ are arithmetically definable, so the set defined by the formula $\phi$ is also arithmetically definable and hence it cannot satisfy the disquotational axioms of $\UTB^-(\Lang_{\PA})$. 
	\end{proof}
	
	In the proof of Theorem  \ref{th_theory_with_rec_sat_without_utb}, we introduced an infinite family of predicates. The ``further'' predicates $T_i$ are from satisfying $\UTB^-(\Lang_{\PA})$, the "closer" $T_{\omega}$ gets to satisfying it. In particular, any completion of $U$ will  actually have a predicate which satisfies $\UTB^-(\Lang_{\PA})$ provably in the theory. In fact, this pattern essentially turns out to hold in the full generality which is the key idea of Theorem \ref{th_theory_with_rec_sat_without_utb}. Let us note that the proof is largely inspired by a related work of \cite{casanovas_farre}. 
	
	\begin{proof} [Proof of Theorem \ref{th_recursive_saturation_implies_truth}]
		
		Let $U$ be a theory in a countable language $\Lang$ extending $\PA$. Suppose that for every model $\M\models U$, the arithmetical part of $\M$ is recursively saturated.		
		Let us introduce some notation. By $\tau(x,y)$ we mean the type consisting of all formulae of the form:
		\begin{displaymath}
			\forall t < x \ \Big( \phi(\val{t}) \equiv \qcr{\phi(t)} \in y \Big),
		\end{displaymath}
		where $\phi$ is an arithmetical formula. In other words, $\tau(x,y)$ expresses "$y$ is a code of the truth predicate for standard formulae and terms not greater than $x$." Observe that for every $\M\models \PA$ and every $a\in M$, $\tau(a,y)$ is a type over $\M$. Moreover, if $\M$ is any model of the elementary arithmetic $\EA$ and for every $a\in M$, $\tau(a,y)$ is realized in $\M$, then $\M\models \PA$. 
		
		Let $\Lang^*$ extend $\Lang$ with one fresh constant $a$. We will inductively build $\Lang^*$-theories $U_{\alpha}$ extending $U$ and sets of $\Lang^*$-formulae $A_{\alpha}$. We let $U = U_0$. For any $\alpha$, let
		\begin{displaymath}
			A_{\alpha} = \set{\exists y \phi(x,y) \in \form_{\Lang}}{ U_{\alpha} \vdash \forall y \left( \phi(a,y) \rightarrow  \tau(a,y) \right)}.
		\end{displaymath}
		The occurrence of $\tau$ should be understood schematically: $U_{\alpha}$ proves every sentence resulting from the above template, where we substitute some formula $\pi \in \tau$ for $\tau$. For any $\alpha$, we let:
		\begin{displaymath}
			U_{\alpha+1} = U_{\alpha} \cup \set{\neg \psi(a)}{\psi(x) \in A_{\alpha}}. 
		\end{displaymath} 
		Finally, we define $U_{\gamma}$ as the union of $U_{\beta}$ for $\beta<\gamma$ for limit ordinals $\gamma$.
		
		Observe that if $U_{\alpha}$ is consistent, then $U_{\alpha+1}$ strictly extends $U_{\alpha}$. Indeed, since  $U_{\alpha}$ extends $U$ and therefore in each model of $U$, $\tau(a,y)$ is realised, this follows immediately  by the Omitting Types Theorem (see \cite[Theorem 2.2.9]{chang-keisler}).

		Now, since $U_{\alpha}$'s form an increasing chain of sets, it has to stabilise. Let $\alpha$ be the least ordinal such that $U_{\alpha} = U_{\alpha+1}$. This means that $U_{\alpha}$ is inconsistent. In particular, $\alpha$ has to be a successor ordinal, say $\alpha = \beta+1$. Since $U_{\beta+1}$ is inconsistent, by compactness there has to be a finite collection $\psi_1(a), \ldots, \psi_n(a) \in A_{\beta}$ such that 
		\begin{displaymath}
			U_{\beta} \vdash \bigvee_{i \leq n} \psi_i(a).
		\end{displaymath}
		
		We now turn to the main part of the proof. Fix any model $\M \models U$. Let $\M^* \succeq \M$ be a countably saturated elementary extension of $\M$.  We claim that for some $\gamma$ and some $\psi \in A_{\gamma}$, $\psi(x)$ holds cofinally in $\M^*$. Indeed, suppose for contradiction that for any ordinal $\delta$, for any $\psi \in A_{\delta}$, there exists $b$ such that $\neg \psi(a)$ holds for all $a>b$. By saturation, this means that for any $\delta$, there exists $b_{\delta}$ such that for all $a>b_{\delta}$, $\neg \psi(a)$ holds for all $\psi \in A_{\delta}$. By  transfinite induction, this means that for all $\delta$, and all $a>b_{\delta}$, $(\M^*, a)\models U_{\delta+1}$. However, as we already noticed, there exists $\beta$ such that $U_{\beta} \vdash \bigvee_{i \leq n} \psi_i(a)$ for a finite family of formulae in $A_\beta$. Since there are finitely many of them, one of these formulae has to hold cofinally in $\M^*$. This is a contradiction.
		
		So fix the minimal $\gamma$ such that one of the formulae in $A_{\gamma}$ holds cofinally in $\M^*$. By minimality, for every $\delta < \gamma$ and every formula $\xi \in A_{\delta}$, the elements satisfying $\xi$ are bounded in $\M^*$. Again by saturation, there exists $b_{\gamma}$ such that no element $a>b_{\gamma}$ satisfies any formula $\xi \in A_{\delta}, \delta < \gamma$. In other words, for every element $a$ greater than $b_{\delta}$, $(\M,a)\models U_{\gamma}$. Fix any such $b_{\gamma}$.
		
		Let $\psi=\exists y\phi(x,y) \in A_{\gamma}$ be an arbitrary formula which is satisfied for cofinally many elements in $\M^*$. Then the following formula $T^*(\xi)$ defines a predicate satisfying $\UTB^-(\Lang_{\PA})$ in $\M^*$:
		\begin{displaymath}
			T^*(\xi):= \forall x \exists y> x\exists z \ \Big( \phi(y,z) \wedge \xi \in z \Big).
		\end{displaymath}
		I.e., $T^*$ says that there are arbitrarily large elements $y$ such that for some $z$, $\phi(y,z)$ holds and $\xi \in z$. We check that $T^*$ satisfies $\UTB^-(\Lang_{\PA})$ in $\M^*$. We work in $\M^*$. Assume first $T^*(\qcr{\eta(t)}))$ holds. Let $a\in M^*$ be greater than $\max\{t, b_{\gamma}\}$. Then there are $b,c>a$ such that $\phi(b,c)$ and $\xi \in c$. Since $\exists y \phi(x,y)\in A_{\gamma}$ and $(\M^*,b)\models U_{\gamma}$, $c$ realizes the type $\tau(b,c)$. It follows by the definition of $\tau(x,y)$ that $\eta(\val{t})$ holds in $\M^*$. Assume conversely that $\eta(\val{t})$ holds. We show that $T^*(\eta(t))$ holds as well. Fix any $a$ and let $b$ be any element greater than $\max\{b_{\gamma}, a, t\}$ such that for some $c$, $\phi(b,c)$ holds. Since $(\M^*, b)\models U_{\gamma}$, it follows that $c$ realizes $\tau(b,y)$ and in particular, $\qcr{\eta(t)}\in c$. Since the initially chosen $a$ was arbitrary, it follows that $T^*(\qcr{\eta(t)})$.
		
		Since $T^*$ is in fact an $\Lang$-formula, by elementarity, $T^*$ witnesses that $\ElDiag(\M)$ syntactically defines $\UTB^-(\Lang_{\PA})$ modulo $\Lang_{\PA}$.
	\end{proof}
	

	Using $\UTB^-$ as a bootstrap we can in fact give an example of the least theory, in the sense of semantic definability, which imposes recursive saturation. Recall the theory $\RSAT(\Lang_{\PA})$ from Remark \ref{rem_theory_of_realisers}.

    In the proof of Theorem \ref{tw_universality_rsat} we shall make use of the following result(\cite[Theorem 4.1]{wcislyk_models_weak})
    \begin{tw}\label{tw_CT_defines_UTB}
    Suppose that $\M\models \PA + \CT^-(a)$, where $a$ is nonstandard. Then there is a formula $\theta(x)\in\Lang_{\M}$ and a nonstandard $b$ such that $\M\models \CT(b)[\theta(x)/T(x)].$
    \end{tw}
	
	\begin{tw}\label{tw_universality_rsat}
		For an arbitrary theory $U\supseteq \PA$ in a countable language $\Lang$ extending $\PA$ the following conditions are equivalent:
		\begin{enumerate}
			\item $U$ imposes $\Lang_{\PA}$-recursive saturation.
			\item $U$ semantically defines $\RSAT(\Lang_{\PA})$ modulo $\Lang_{\PA}$.
		\end{enumerate}
	\end{tw}
	\begin{proof}
		The implication $2.\Rightarrow 1.$ is straightforward, so we prove the reverse direction. Assume $U$ imposes recursive saturation and take any $\M\models U$. Let $T(x)$ be a definable $\UTB^-(\Lang_{\PA})$-predicate and let $\Lang_T$ denote the language $\Lang_{\PA}\cup \was{T}$. For an arbitrary primitive recursive type $p(x,y)$, let $\was{\phi_i(x,y)}_{i\in\omega}$ be its p.r. enumeration. To enhance readability we will treat the function $i\mapsto \phi_i$ as a new primitive symbol (that should be eliminated by substituting the proper absolute $\Sigma_1$ representation of this function). Let Prog$_z$[$P(z,\bar{x})$] be the following formula ($P(z,\bar{x})$ is a placeholder for an arbitrary formula):
		\[\forall z \bigl(\forall i<z P(i,\bar{x})\rightarrow P(z,\bar{x})\bigr).\]
		We distinguish two cases:
		\paragraph{Case 1: The standard cut is definable in $\M$.} 
		
		Let $\omega(x)$ be a $\Lang(\M)$ definition of the standard cut of $\M$. 
		Our definition of $R_p(x,y)$ formalizes the following definition by cases:
		\begin{displaymath}
				\begin{cases}\forall z \bigl(\omega(z)\rightarrow T(\phi_z(\dot{x},\dot{y})\bigr), &\textnormal{if } \textnormal{Prog}_v[\exists x\forall z<vT(\phi_z(\dot{x},\dot{y}))] \\ 
					\forall z<vT(\phi_z(\dot{x},\dot{y})), &\textnormal{if $v$ is the greatest such that }\bigcap_{z<v}\set{x}{T(\phi_z(\dot{x},\dot{y}))}\neq \emptyset
				\end{cases}
		\end{displaymath} 
		We now prove that this definition works. That the optimality axioms hold is clear: fix any parameter $y$ and any $n\in\omega$ and suppose that $\bigcap_{z<n}\set{x}{T(\phi_z(\bar{x},\bar{y}))}\neq \emptyset$. Then regardless of which case of our definition holds, any $x$ satisfying $R_p(x,y)$ belongs to the set $\bigcap_{z<n}\set{x}{T(\phi_z(\dot{x},\dot{y})}$, which, by the axioms of $\UTB^-(\Lang_{\PA})$, is the same as $\set{x}{\bigwedge_{z<n}\phi_z(x,y)}$. Let us deal with the  non-emptiness axioms. Fix any $y$. If $\neg \textnormal{Prog}_v[\exists x\forall z<vT(\phi_z(\dot{x},\dot{y}))]$, then by definition we have $\exists x R_p(x,y)$. Suppose now $\textnormal{Prog}_v[\exists x\forall z<vT(\phi_z(\dot{x},\dot{y}))]$ holds.
		Then $p(x,y)$ is realised, since $\M$ is recursively saturated. In particular $\exists x \forall z \bigl(\omega(z)\rightarrow T(\phi_z(\dot{x},\dot{y})\bigr)$, which concludes the (current case of the) proof.
		
		\paragraph{Case 2: The standard cut is not definable in $\M$.} 
		
		By the compositional axioms of $\UTB^-$ we have that for every $n\in\omega$ $(\M,T^{\M})\models\CT^-(n)$. So, by our assumption there is a nonstandard number $a$ such that $(\M, T^{\M})\models \CT^-(a)$, i.e. $T$ is a compositional truth predicate for all sentences of depth at most $a$. Then, by Theorem \ref{tw_CT_defines_UTB}, 
		there is a $\Lang_{M}$- formula $T'(x)$ which defines a partial inductive satisfaction class, i.e. for some nonstandard $b$
		\[(\M, (T')^{\M})\models \CT^-(b) + \Ind_{\Lang_{T}}.\]
		Let us fix $T'$ and $b$ such that the above holds. Let $c$ be any nonstandard number such that $\M\models \forall i<c \ \phi_i\in\dpt(b)$. Now we can simply define $R_p(x,y)$ to be the following formula:
		\rem{\begin{equation*}
				\exists v \bigl[\forall i <v T'(\phi_i(\dot{x},\dot{y})) \wedge  \forall v'>v \bigl(v'>c \vee \neg \exists z\forall i<v' T'(\phi_i(\dot{z},\dot{y}))\bigr) \bigr].
		\end{equation*}}
		Thus $R_p(x,y)$ defines the set $\bigcap_{z<v}\set{x}{T'(\phi_z(\dot{x},\dot{y}))}$, where $v$ is the greatest number $\leq c$ such that  the above intersection is not empty (by convention the intersection of the empty family is the whole universe of the model). Such a $v$ always exists by the fact that $(\M, (T')^{\M})\models \Ind_{\Lang_{T}}$. It is now a routine exercise to check that for this definition both the nonemptiness and optimality axioms hold.
	\end{proof}

	\begin{uwaga}
		It is easy to observe that Theorems \ref{th_recursive_saturation_implies_truth} and \ref{tw_universality_rsat} still hold if $\PA$ and $\Lang_{\PA}$ are substituted by a pair $V$, $\Lang$ such that $\PA\subseteq V$, $\Lang$ has at most countably many additional predicates and $V\vdash \Ind_{\Lang}$.
	\end{uwaga}

	By modifying the proof of Pakhomov and Visser (\cite[Theorem 4.1]{pakhomov_visser}) we can show that actually no finite theory can replace RSAT in Theorem \ref{tw_universality_rsat}:
	
	\begin{stwierdzenie}
		If $U$ is a finite theory such that $\PA+U$ imposes $\Lang_{\PA}$-recursive saturation, then there is a finite theory $V$ such that $\PA + V$ does not semantically define $U$ modulo $\Lang_{\PA}$ and $\PA+V$ imposes $\Lang_{\PA}$-recursive saturation.
	\end{stwierdzenie}
	\begin{proof}
		In the proof all definability between theories is modulo $\Lang_{\PA}$. First observe that for a finite theory $U$ and an arbitrary theory $V$, $V$ semantically defines $U$ if and only if $V$ syntactically defines $U$. Fix $U$ and assume that $U$ imposes $\Lang_{\PA}$-recursive saturation. For an arithmetical formula $\phi(x)$ let $\CT^-\restr{\phi(x)}$ denote the following $\Lang_T$ sentence:
		\[\forall x \bigl(\forall y<x \ \phi(y)\rightarrow \CT^-(\Lang_{\PA})(x)\bigr).\]
		Now apply the diagonal lemma to get a $\Sigma_1$ formula $\psi = \exists x \phi(x)$ such that
		\[\EA\vdash \psi\equiv "\PA + \CT^-\restr{\neg\phi(x)}\textnormal{ syntactically defines } U".\]
		We shall show that $\CT^-\restr{\neg\phi(x)}$ satisfies the requirements put on a theory $V$ in the thesis of our proposition. First, we claim that $\psi$ is false. Assume the opposite and take $n\in\omega$ such that $\EA\vdash \phi(\num{n})$. By the existence of partial truth predicates, it follows that $\EA$ defines $\EA+ \CT^-\restr{\neg\phi(x)}$. However, by the fixpoint property of $\psi$, it follows that $\PA$ defines $U$. Since $U+\PA$ imposes $\Lang_{\PA}$-recursive saturation, it follows that $\PA$ imposes $\Lang_{\PA}$-recursive saturation, which is obviously false.
		Hence $\psi$ is false and it follows that $\PA+\CT^-\restr{\neg\phi(x)}$ does not define $U$. We claim that $\PA+\CT^-\restr{\neg\phi(x)}$ imposes $\Lang_{\PA}$ recursive saturation. Indeed, fix any $\M\models \PA + \CT^-\restr{\neg\phi(x)}$. Since $\psi$ is false, then for every $n\in\omega$, $\M\models \neg\phi(\num{n})$, hence by the overspill principle, for some nonstandard $c$, $\M\models \forall x<c \neg\phi(x)$. As a consequence, $\M\models \CT^-(c)$. It follows by Lachlan's theorem (see \cite[Theorem 15.5]{kaye}) that $\M\restr{\Lang_{\PA}}$ is $\Lang_{\PA}$-recursively saturated, which ends the proof.
	\end{proof}
	
	Finally, we would like to return to Theorem \ref{th_theory_with_rec_sat_without_utb}. We will show that the infinite language we used in the proof of that theorem was not an indispensable ingredient of our argument. In the refined proof, instead of saying that an alleged definition of a truth predicate cannot use all the symbols, we will only say that it cannot have a bounded complexity. Thus, this shall fullfil our promise made in the introduction of delivering a theory in a finite language which defines $\UTB^-(\Lang_{\PA})$ semantically but not syntactically. 
	
	In the proof, we will use properties of generic sets in computability theory. Our argument consists of two parts. We will first show in Lemma \ref{lem_hard_coding_in_N} that for an arbitrary $n$, we can find a set $A \subseteq \mathbb{N}$ (in the standard model) such that it defines the standard truth predicate for the arithmetical sentences, but not with a formula of complexity $\leq n$. To accomplish this, we amalgamate two constructions: (iterated) jump-inversion theorem and the construction of an $n$-generic (as presented in \cite[Chapter V]{odifreddiI} and \cite[Chapter XII]{odifreddiII}). Then we will show that this argument formalises in arithmetic and we will use Lemma \ref{lem_hard_coding_in_N} to construct a theory in a finite language which defines a $\UTB^-(\Lang_{\PA})$ predicate in every model but does not do so uniformly. Below $\leq_T$, $\equiv_T$ denote Turing-reducibility and Turing-equivalence, respectively. 
	
	Let us recall the definition and the basic properties of an $n$-generic set.

	\begin{definicja} \label{def_arithmetical_forcing}
		Let $s$ be a finite binary sequence. We define the forcing relation $s \Vdash \phi$ for $\phi$ first order formulae with a single second-order variable $X$ by induction on the complexity of $\phi$:
		\begin{itemize}
			\item If $\phi$ is a first-order atomic formula, then $s \Vdash \phi$ iff $\mathbb{N} \models \phi$.
			\item If $\phi = n \in X$, then $s \Vdash \phi$ iff $s(n) = 1$ (in particular $|s| > n$).
			\item If $\phi = \neg \psi$, then $s \Vdash \phi$ iff for all $t \supseteq s$, $t \nVdash \psi$.
			\item If $\phi = \psi \vee \eta$, then $s \Vdash \phi$ iff $s \Vdash \psi$ or $s \Vdash \eta$. 
			\item If $\phi = \exists x \psi(x)$, then $s \Vdash \phi$ iff for some $n \in \mathbb{N}$, $s \Vdash \psi(n)$. 
		\end{itemize} 
	\end{definicja}

    \begin{definicja}
    Let $A\subseteq \mathbb{N}$.
    \begin{itemize}
        \item We say that $A$ forces $\phi$, $A\Vdash \phi$ iff for some finite $s\subseteq A$, $s\Vdash \phi$.
        \item $A$ is $n$-generic iff for every $\Sigma^0_n$ sentence $\phi$ $A\Vdash \phi$ or $A\Vdash \neg\phi$.
    \end{itemize}
    \end{definicja}

    \begin{tw}
    A set $A\subseteq \mathbb{N}$ is $n$-generic iff for every $\phi\in\Sigma^0_n \cup \Pi_n^0$, $(\mathbb{N},A)\models \phi$ iff $A\Vdash \phi$.
    \end{tw}
	
	\begin{lemat} \label{lem_hard_coding_in_N}
		For any $n \in \omega$, there exists a set $A \subseteq \mathbb{N}$ such that $\Th(\mathbb{N})\subseteq\mathbb{N}$ is definable in $(\mathbb{N},A)$, but $\Th(\mathbb{N})$ is not definable in $(\mathbb{N},A)$ with a formula of complexity $< \Sigma_n$.
	\end{lemat}
	\begin{proof}
		For a given $n$, we will find a set $A$ such that:
		\begin{itemize}
			\item[a.] $A \oplus 0^{(n)}$ computes $0^{(\omega)}$.
			\item[b.] $A \oplus 0^{(n)}$ computes $A^{(n)}$. 
		\end{itemize}
		We claim that if $A$ satisfies a. and b., then $\Th(\mathbb{N})$ is not definable in $(\mathbb{N},A)$ with a formula of complexity $k< n$.
		 If it is, then $0^{(\omega)} \leq_T  A^{(k)}.$ This yields a contradiction as follows:
		\begin{displaymath}
			A^{(n)} \leq_T 0^{(n)} \oplus A \leq_T 0^{(\omega)} \oplus A \leq_T  A^{(k)}. 
		\end{displaymath}
		
		Condition a. will be guaranteed by the construction from the Jump Inversion Theorem (see \cite[Theorem V.2.24]{odifreddiI}) while b. will be guaranteed by the construction of an $n$-generic (see \cite[XII.1.11]{odifreddiII}). 
		As stated in \cite[Proposition XII.1.6]{odifreddiII}, the restriction of the forcing relation to $\Sigma_n$ sentences is computable in $0^{(n)}$. We fix an enumeration $\phi_i(X)$ of $\Sigma^0_n$ formulae with one second-order variable 
		 and an enumeration $\xi_i$ of all $\Lang_{\PA}$ sentences.
		
		We construct $A$ in stages. At each stage $n$, we have a finite binary sequence $s_n$. We proceed as follows:
		\begin{itemize}
			\item If $n+1$ is even, say $2k$, then we look at the formula $\phi_k$. If there exists an extension $s \supsetneq s_n$ such that $s\Vdash \phi_k$ holds, then we let $s_{n+1}$ be the lexicographically least such $s$. If not, we let $s_{n+1} :=s_n$. 
			\item If $n+1$ is odd, say $2k+1$, then we set $s_{n+1}$ equal to $s_n\frown \tuple{1}$ if $\xi_k$ is true and $s_n\frown \tuple{0}$ if it is false.
		\end{itemize}
		Finally, we let $s= \bigcup s_n$ and we let $A$ be the only set whose characteristic function is $s$. Clearly, $A$ is $n$-generic.
		
		Indeed, assume that $A\nVdash \phi$ and $\phi = \phi_k$. By definition, for every $s\subseteq A$, $s\nVdash \phi_k$. In particular $s_{2k}\nVdash \phi_k$, hence $s_{2k} = s_{2k-1}$. Moreover no $\tau\supset s_{2k}$ forces $\phi_k$, so in particular we get that $s_{2k}\Vdash \neg\phi_k$. This concludes the proof of $n$-genericity of $A$. We want to check that $A$ satisfies both conditions a. and b. 
		
		\paragraph*{Claim I: $A \oplus 0^{(n)}$ computes $0^{(\omega)}$.}
		
		The set $0^{(\omega)}$, the simple join of $0^{(k)}$ for $k < \omega$ is the Turing degree of the standard arithmetical truth predicate $\Th(\mathbb{N})$. We will show how to compute the latter set from $A \oplus 0^{(n)}$. Fix any $k$. We want to decide whether $\xi_k \in \Th(\mathbb{N})$. We inductively compute the sequence of finite sequences
		\begin{displaymath}
			s_0 \subseteq s_1 \subseteq \ldots \subseteq s_{2k+1}.
		\end{displaymath}
		Suppose that we know $s_{2m}$. Then we compute $s_{2m+1}$ as the only one-step extension of $s_{2m}$ agreeing with $A$. On the other hand, if we know $s_{2m-1}$, we can define $s_{2m}$ as the lexicographically smallest sequence extending $s_{2m-1}$ and forcing $\phi_m$. 
		
		Finally, we set $\xi_k \in \Th(\mathbb{N})$ iff $s_{2k+1}$ extends $s_{2k}$ with a single digit $1$. 
		
		\paragraph*{Claim II: $A^{(n)} \equiv_T A \oplus 0^{(n)}$.}
		
		It is enough to show that $A^{(n)} \leq_T A \oplus 0^{(n)}$. We want to check, for any $\Sigma_n^0(A)$ formula $\phi$ whether $\phi(A)$ holds. As noted before, we can compute the sequence $s_0, s_1, s_2, \ldots$ from $A\oplus 0^{(n)}$. Assume that $\phi = \phi_k$. We claim that $\phi_k(A)$ holds iff $s_{2k}$ nontrivially extends $s_{2k-1}$. Since $A$ is $n$-generic, we can reason as follows: 
		\begin{align*}
			(\mathbb{N},A)\models \phi_k(A) &\equiv A\Vdash \phi_k\\
			&\equiv \exists s\subset A\ \ s\Vdash \phi_k\\
			&\equiv s_{2k}\Vdash \phi_k\\
			&\equiv s_{2k-1}\subsetneq s_{2k}
		\end{align*}
		The first equivalence is by \cite[Proposition XII.1.10]{odifreddiII}, the second is a definition, the third and fourth follow by the construction. Now, notice that the sequences $s_i$ can be computed using  $A \oplus 0^{(n)}$ as an oracle. Indeed, in order to compute $s_{2k}$ from $s_{2k-1}$, we only need to know whether there exists $s \supsetneq s_{2k-1}$ such that $s \Vdash \phi_k$ which, as we already mentioned, is a $\Sigma_n$ fact. In order to compute $s_{2k+1}$ from $s_{2k}$, we only have to check whether the corresponding single digit of the set $A$ is $0$ or $1$. This concludes the argument and completes the proof of the lemma. 
	\end{proof}
	
	In Lemma \ref{lem_hard_coding_in_N}, we produced some subsets of $\mathbb{N}$ which define $\Th(\mathbb{N})$ but only with a formula no simpler than a fixed complexity $k$. We will check that the construction can actually be captured by a \emph{theory}, which allows us to improve Theorem \ref{th_theory_with_rec_sat_without_utb} to theories in finite languages. 
	
	\begin{tw} \label{th_theory_finlang_with_rec_sat_without_utb}
		There exists a theory $U$ in a finite language extending $\PA$ which imposes $\Lang_{\PA}$-recursive saturation but does not syntactically define $\UTB^-(\Lang_{\PA})$ modulo $\Lang_{\PA}$.
	\end{tw}
	
	\begin{proof}
		As previously, in the proof all definability between theories is modulo $\Lang_{\PA}$. Let us fix a primitive recursive enumerations $(\phi_i)$ of arithmetical formulae and, for each $k$, $(\psi^k_j)$ of arithmetical $\Sigma_k$ sentences with a single second order variable $X$. We shall construct a theory $U$ in the arithmetical language expanded with two predicates $A$, $T_\omega$, and a constant $c$. The axioms of $U$ will comprise the axioms of $\PA$, full scheme of induction for the extended language and an infinite list of axioms described below.
		
		We first define an auxiliary sequence of theories. For each $k$, the theory $\alpha_k$ says that $A$ arises from a construction of the $k$-generic code of $\Th(\mathbb{N})$. More precisely, $\alpha_k$ consists of sentences $\alpha_k^{\phi(x)}$, for $\phi(x)\in\Lang_{\PA}$, where $\alpha_k^{\phi(x)}$ says: 
		\begin{center}
			,,For each $x$, there exists an $m$ and a sequence of finite sequences $\emptyset = s_0 \subseteq s_1 \subseteq \ldots \subseteq s_m$ satisfying the conditions specified below, such that for all $y \leq x$, $A(y)$ iff $s_m(y) =1$, where the conditions on the sequences $s_{i+1}$ are as follows:
			\begin{itemize}
				\item Case I: if $j = 2l$ for some $l$, then
				\begin{itemize}
					\item If there exists a sequence $s \supsetneq s_{j}$ such that $s\Vdash \psi^k_j$, then $s_{j+1}$ is the lexicographically minimal such $s$ strictly extending $s_j$.
					\item $s_{j+1} = s_j$, otherwise. 
				\end{itemize}
				\item Case II: if $j = 2l+1$, then $s_j = s_{2l} \frown \tuple{i}$ for $i \in \{0,1\}$; moreover, if there exists a closed term $t$ such that $\phi_j = \phi(t)$,  then:
				\begin{itemize}
					\item If $\phi(\val{t})$, then $s_{j+1} = s_j \frown \tuple{1}$.
					\item If $\neg \phi(\val{t})$, then $s_{j+1} = s_j \frown \tuple{0}$.
				\end{itemize} 
			\end{itemize}
		\end{center}
		Note that, while in Case I there is no bound on the number of bits $s_{j+1}$ adds to $s_j$, in Case II $s_{j+1}$ is always an extension of $s_j$ by a single bit.
		
		Finally, $U$ is the union of the following sentences:	
		\begin{itemize}
			\item $	\PA$.
			\item The full induction scheme for the extended language.
			\item $c=\num{k} \rightarrow \alpha_k^{\phi(\bar{x})}$, for $\phi$ formulae of $\Lang_{\PA}$.
			\item $	c > \num{3m} \rightarrow \forall t\Big(T_{\omega}(\qcr{\phi_m(t)}) \equiv \phi_m(\val{t})\Big)$.
		\end{itemize}
		
		It is enough to check that $U$ defines $\UTB$ semantically but not syntactically.
		
		\textbf{Claim I: $U$ defines $\UTB(\Lang_{\PA})$ semantically}
		
		Let $\M \models U$. Then, in $M$, $c$ is either a standard or a nonstandard number. If $c$ is nonstandard, then $T_{\omega}$ satisfies $\UTB$. If $c = k \in \omega$, then fix any standard arithmetical formula $\phi$ and a closed term $t\in M$ (perhaps nonstandard). We desrcibe an $\M$-definable procedure for checking whether $\phi(\val{t})$ holds in $\M$.
		
		Fix $d$ such that $\phi(t) = \phi_d$ in our enumeration. We perform the construction of $A$. Notice that initial segments of the sequence $s_0 \subseteq s_1 \subset \ldots$ defined in the axiom $\alpha^{\phi(x)}_k$ are actually definable from $A$: the even steps are computable in the sense of $\M$ using as the oracle the definable set of all true $\Sigma_k$ sentences (with parametres). The odd steps are directly definable in the oracle $A$ (we simply check whether $A$ has $0$ or $1$ in the respective places). In particular, we set $T(\phi(t))$ iff the last bit of $s_{2d+1}$ is $1$. 
		
		\textbf{Claim II: $U$ does not syntactically define $\UTB(\Lang_{\PA})$}
		
		Suppose that $U$ syntactically defines $\UTB(\Lang_{\PA})$ with a formula $\phi$. Suppose that $\phi \in \Sigma_{k}$. Consider the model $\mathcal{N}$ whose arithmetical part is the standard model $\mathbb{N}$, the constant $c$ is interpreted as $3k$, $T_{\omega}$ is interpreted as the partial arithmetical truth predicate for the arithmetical $\Sigma_{k}$ formulae (which can be assumed to include $\phi_0, \ldots, \phi_{k-1}$) and $A$ is the set constructed in the proof of Lemma \ref{lem_hard_coding_in_N} for $n= 3k$.	Then obviously $\mathcal{N}\models U$. 
		
		In the $\Sigma_k$ formula $\phi$, we can substitute the partial arithmetical predicate $\Tr_{\Sigma_{k}}$ for $T_{\omega}$ and $\num{3k}$ for $c$, obtaining a $\Sigma_{2k+1}$ formula $\phi^*$ which defines $\Th(\mathbb{N})$ in $(\mathbb{N}, A)$. However, by Lemma \ref{lem_hard_coding_in_N}, for any $l<3k$, $\Th(\mathbb{N})$ cannot be defined with a $\Sigma_{l}$ formula from the predicate $A$. 
	\end{proof}

	\section{Summary and open problems}
	
	As we mentioned in the introduction, the questions which we tackle in this article have several different ``dimensions''. In order to facilitate navigating between them, we summarise our findings below.
	
	The table below concerns the problem which theories impose truth-like semantic properties. 
	\begin{center}
		
		\begin{tabular}{c||c|c||c|c||c|c||c}
			
			&\multicolumn{6}{c||}{Sequential theories, finite languages} & $\PA$, ctble languages \\ 
			\cmidrule{2-8}
			& \multicolumn{2}{c||}{ImpElemEq} & \multicolumn{2}{c||}{PreservDef} & \multicolumn{2}{c||}{ImpElem} & \multirow{2}{*}{\cellcolor{gray!10}ImpRecSat}\\ \cmidrule{2-7}
			& \cellcolor{gray!10}non-Uni & \cellcolor{gray!40}Uni & \cellcolor{gray!10}non-Uni & \cellcolor{gray!40}Uni &\cellcolor{gray!10}non-Uni & \cellcolor{gray!40}Uni &\cellcolor{gray!10} \\ \toprule
			Sem. defi. & \cellcolor{gray!10}$\TB$ & \cellcolor{gray!40}$\TB$ & \cellcolor{gray!10}$\DEF$  &  \cellcolor{gray!40}$\DEF$ & \cellcolor{gray!10}$\USB$ & \cellcolor{gray!40}$\USB$ & \cellcolor{gray!10}$\textnormal{RSAT}$\\ \hline
			Synt. defi. & none & \cellcolor{gray!40}$\TB$ & none & \cellcolor{gray!40}$\DEF$ & none & \cellcolor{gray!40}$\USB$ & none	\\
			\bottomrule
		\end{tabular}
	\end{center}
	The first row contains information	for which theories our classification works. The second row contains (the abbreviations of) the main model-theoretic properties considered in this paper. The first three come in two variants: uniform and non-uniform. The row "Sem. (Synt.) defi." contains theories that are characterized by the respective property up to semantic (syntactic) definability. When we write "none," this means that none \textit{of the main theories} is characterized by the respective property. In all the cases the theory defined in the proof of Theorem \ref{th_theory_finlang_with_rec_sat_without_utb} serves as a counterexample. Moreover, in the semantic sense and for extensions of $\PA$, the non-uniform versions of the properties get properly more restrictive from left to right (light grey properties and theories). For instance, every theory which imposes $\Lang_{\PA}$-recursive saturation (for extensions of $\PA$) preserves $\Lang_{\PA}$-definability, but there is a theory ($\DEF$) which preserves $\Lang_{\PA}$-definability but does not impose $\Lang_{\PA}$-recursive saturation. The same holds over all natural sequential theories, except for some degenerate cases, for uniform versions of model-theoretical properties and their respective theories (dark grey properties and theories).

	In this article, we have also analysed and investigated relations between the axiomatic theories corresponding to the semantic properties which truth-like predicates can impose. We would now like to lay down a diagram depicting the interdependencies between those theories. 
	
	Below, by writing a filled-arrow from $U$ to $V$, we mean that $V$ defines $U$ \textbf{syntactically} (think of that as a statement that $U$ embeds into $V$). Whenever we \emph{do not} write an arrow from $V$ to $U$ this means that $V$ does not define $U$ \textbf{semantically}, over any sequential theory. The only exception is when there is a path consisting of other arrows. If we \emph{do not know} whether $V$ \textbf{semantically} defines $U$ or this happens sometimes but not in typical cases, we write a dashed arrow with a question mark (which always means that we do not even know whether semantic definability obtains) and, possibly, some comments as to what is known.
	
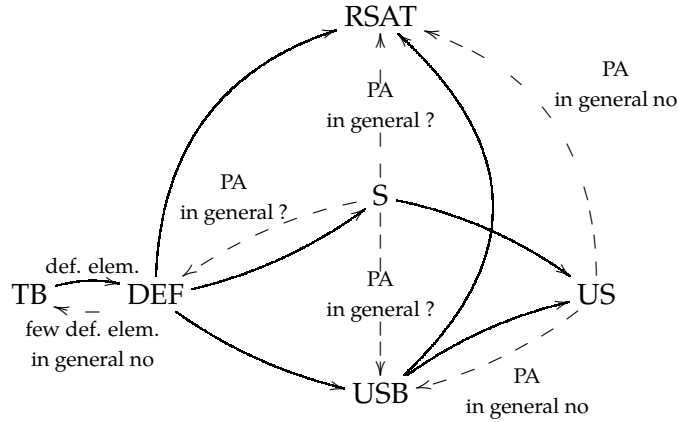
\begin{displaymath}
		\xymatrix{
			&  & & \RSAT  & & \\
			& & & & & \\
			& &  & \textnormal{S} \ar @{-->} [dd]|{\txt{\scriptsize{PA} \\ \scriptsize{in general ?}}} \ar @{-->} [uu]|{\txt{\scriptsize{PA} \\ \scriptsize{in general ?}}} \ar @/^/ [drr]  \ar @{-->} @/_/ [lld]_{\txt{\scriptsize{PA} \\ \scriptsize{in general ?}}} & & \\
\TB \ar @/^/ [r]^{\txt{\scriptsize
		{def. elem.}}} & \DEF \ar @/^2.3pc/ [uuurr] \ar @/^/ @{-->} [l]^{\txt{\scriptsize{few def. elem.} \\ \scriptsize{in general no}}}	\ar @/_/ [urr] \ar @/_/ [drr]	& & & & \textnormal{US} \ar @{-->} @/_2.3pc/ [uuull]_{\txt{\scriptsize{PA} \\ \scriptsize{in general no}}} \ar @{-->} @/^/ [lld]^{\txt{\scriptsize{PA} \\ \scriptsize{in general no}}} \\
 	& & & \USB \ar@/_3.5pc/  [uuuu]  \ar @/^/ [urr] & &
		}
	\end{displaymath}
Most of the entries in the diagram, both positive and negative are directly stated in this article, mostly in Section \ref{sec_definability_between_theories}. The fact that there is no arrow from $\textnormal{US}$ to $\textnormal{S}$ follows from the fact that $\textnormal{US}$ naturally defines $\textnormal{USB}$ and by Proposition \ref{stw_S_nie_definiuje_USB}, we know that over no sequential theory can $\textnormal{S}$ define $\USB$ semantically, so any arrow from $\textnormal{US}$ to $\textnormal{S}$ would yield a nonexistent arrow between $\USB$ and $\textnormal{S}$.

The theory $\RSAT$ is a bit of an outlier.\footnote{In this article, we have defined $\RSAT$ only in the context of $\PA$. However, one could formulate a definition suitable for arbitrary sequential theories, following the cases of other theories, like $\USB$.} First, notice that over an arbitrary sequential theory $\RSAT$ imposes recursive saturation which in general none of the other listed theories does, so it is not semantically definable over any of them. 

On the other hand, it can be easily checked that $\RSAT_K$ syntactically defines $\USB_K$ (hence also $\DEF_K$). Namely $S(\qcr{\phi(x)},y)$ holds iff $\forall x \bigl(R(p_{\phi(y)},x,y)\rightarrow \num{1}(x)\bigr)$, where $p_{\phi(y)}$ is the (easily definable from $\qcr{\phi(y)}$)) one-element type
\begin{displaymath}
	(\phi(y)\wedge \num{1}(x)) \vee (\neg\phi(y)\wedge \num{0}(x)).
\end{displaymath}
It follows that, over theories containing $\PA$, $\RSAT_{id}$ defines $\textnormal{US}_{id}$. Moreover, the proof of Theorem \ref{tw_nondef_US_in_UTB} can be adapted to show that over theories in $\Lang_{\PA}$ consistent with $B\Sigma_n+\neg I\Sigma_n + \exp$, $\RSAT_{id}$ does not semantically define $\textnormal{US}_{id}$: given a recursively saturated model $\M\models B\Sigma_n+\neg I\Sigma_n + \exp$ one finds a $A$ set which is definable with a parameter $a\in M$ but which does not contain any element definable from $a$. Then one canonically expands $\M$ to a model of $\RSAT_{id}$ by interpreting $R(p,x,b)$ to be the set of all elements which realize the maximal fragment of $\was{\phi_i^p(x,b)}$. The automorphism group of the expanded model is the same as that of $\M$ and one can simply rewrite the rest of the proof. In the case of general theories, essentially nothing else is known.

Below we list some open questions and perspectives for further research.
	\begin{enumerate}
		\item[Q1] In this article, we considered the notion of definability preservation (so arithmetically definable elements satisfy the same definitions between submodels of a given theory). We also considered theories which impose equality of definables, but we were unable to provide a good characterisation of such theories. One possible candidate is the theory $\DEF$ of definability. 
		
	\textit{
			Let $U$ be any sequential theory in a language $\Lang$ and $K: S^1_2\lhd U$. Is there an r.e. sequential theory in a finite language which imposes equality of $\Lang$ definables but does not semantically define $\DEF_K(\Lang)$ modulo $\Lang$? If the answer is negative, how about syntactic definability of $\DEF_K(\Lang)$ (modulo $\Lang$)?
	}

	\item[Q2] Most questions about relationships between semantic properties and defining respective theories investigated in this paper found a definite answer. However, this is not quite so in the case of the recursive saturation. There are three reasons for that: first of all, we assume that our theory contains full $\PA$, and consequently we loose a lot o generality. Secondly we are only able to obtain semantic, rather than syntactic definability of $\RSAT$ which makes our characterisation deficient in a sense, because it is not a nice ``completeness-style'' result of the kind which holds for other semantic notions of interest. However, our counterexamples seem to make the essential use of the fact that the theory in question is infinite.  Thirdly, our characterisation works only for theories in countable languages (because of the use of the Omitting Types Theorem). These three issues lead to the following questions: 
	
	\textit{Assume that $U$ is a finite theory with $\PA\subset U$ which imposes $\Lang_{U}$-recursive saturation. Does $U$ define $\RSAT_{id}(\Lang_{\PA}$) syntactically?}
	
	\textit{Assume that $U$ is a sequential theory in a countable language $\Lang_U \supset \Lang$ with $K: S^1_2 \lhd U$ which imposes $\Lang$-recursive saturation. Does $U$ define $\RSAT_{K}(\Lang)$ semantically? Assume that $U$ is finite. Does it define $\RSAT_{K}(\Lang)$ syntactically?}
	
	\textit{Assume that $U$ is an arbitrary theory with $\PA\subset U$ (not necessarily in a countable language) which imposes $\Lang_{\PA}$-recursive saturation. Does $U$ define $\RSAT_{id}(\Lang_{\PA}$) semantically?}
	
		\item[Q3] We know that for theories $U$ extending $\PA$ with the identity embedding of $S^1_2$, $\USB_{id}(\Lang_{U})$ defines $\textnormal{S}_{id}(\Lang_{U})$. However, the use of $\PA$ seems essential for this definability and we can expect that the methods used to separate $\USB$ from $\textnormal{US}$ can be modified to get this stronger separation. We note that by Theorem \ref{stw_nondef_H_usb} over any finite fragment of $\PA$, $\USB_{id}$ does not \textit{parameter-free} semantically define $S_{id}$. However, at this point we cannot exclude the existence of a definition with parameters.
		
		\textit{Let $U$ be an r.e. sequential theory and let $K:S^1_2\lhd U$. Does $\USB_K(\Lang_U)$ semantically define $\textnormal{S}_K(\Lang_U)$ modulo $\Lang$ (possibly with parametres)? Does it hold if we take $U$ to be a finite fragment of $\PA$ and $K$ to be the identity embedding?} 

		\item[Q4] One arrow from our diagram of the interdependencies between theories is missing and it is really not clear what kind of separation should hold between them. 
		
		\textit{Let $U$ be an r.e. sequential theory and let $K:S^1_2\lhd U$. Does $\DEF_K(\Lang_U)$ semantically define $\textnormal{S}_K(\Lang_U)$ modulo $\Lang_U$?}
		
		\item[Q5] As we have already mentioned, the status of $\RSAT$ in the diagram of definability between theories is completely open, besides rather trivial observations. 
		
		\textit{Let $U$ be an r.e. sequential theory and let $K:S^1_2\lhd U$. Does $\RSAT_K(\Lang_U)$ semantically define $\textnormal{S}_K(\Lang_U)$ modulo $\Lang_U$? Does it define $\textnormal{USB}_K(\Lang_{U})$?}

		\item[Q6] In this article, we have investigated relations between semantic properties of theories and defining truth-like predicates as well as definability relations between those predicates. Those results combined allow us to obtain some implication between various semantic properties themselves. For instance, if $U$ contains $\PA$ and imposes $\Lang_{\PA}$-recursive saturation, then it also imposes $\Lang_{\PA}$-elementarity between the models. While we have not systematically investigated this kind of relations, we observed that at least some of them lead to nontrivial questions, especially if we impose additional conditions, such as finiteness of theories. Possibly more questions in this spirit can be formulated.

		 \textit{Let $\Lang$ be a language. Is there a finite sequential theory which (uniformly) imposes $\Lang$-elementary equivalence but does not (uniformly) preserve $\Lang$-definability? Same question for the pair preserving definability/imposing elementarity.}
	
		\item[Q7] Above, we have already asked whether $\RSAT$ can be defined semantically in general sequential theories imposing recursive saturation. The following question approaches this problem from a somewhat different angle.
		
		\textit{Let $\Lang$ be a language. Is there a finite sequential theory which imposes $\Lang$-recursive saturation but does not uniformly impose $\Lang$-elementarity?}
	
		\item[Q8] In Subsection \ref{subsect_PA_def_inno}, we analysed definability between the predicates satisfying theories of truth-like notions in specific models. We established that in every nonstandard model of $\PA$, there exists a predicate satisfying $\TB^-$ which does not allow us to define a predicate satisfying $\USB^-$ and, similarly, that predicates satisfying $\DEF^-$ are definable from predicates satisfying $\TB^-$ only in the models of $\Th(\mathbb{N})$. This sort of questions about interdefinability of truth-like theories can be extended to other theories we introduced in this paper. For instance, we can ask about the relationship between $\DEF^-$ and $\USB^-$.
		
		\textit{For which models $\M\models \PA$ it holds that for every $D\subseteq M^2$ such that $(\M,D)\models \DEF^-(\Lang_{\PA})$, $\ElDiag(\M,D)$ defines $\UTB^-(\Lang_{\PA})$? We conjecture that these are precisely the prime models of $\PA$.}
	\end{enumerate}

\end{document}